\documentclass{amsart}

\usepackage{amssymb}
\usepackage{amscd}

\newtheorem{theorem}{Theorem}[section]
\newtheorem{lemma}[theorem]{Lemma}
\newtheorem{lemma-definition}[theorem]{Lemma-Definition}
\newtheorem{proposition}[theorem]{Proposition}
\newtheorem{corollary}[theorem]{Corollary}

\theoremstyle{definition}
\newtheorem{definition}[theorem]{Definition}
\newtheorem{remark}[theorem]{Remark}
\newtheorem{notation}[theorem]{Notation}
\newtheorem{example}[theorem]{Example}

\numberwithin{equation}{section}

\newcommand{\co}{\mathop{\mathrm{co}}}
\newcommand{\gr}{\mathop{\mathrm{gr}}}
\newcommand{\ch}{\mathop{\mathrm{ch}}}
\newcommand{\Cleft}{\mathrm{Cleft}}
\newcommand{\ZCleft}{\mathrm{ZCleft}}
\newcommand{\Ker}{\mathop{\mathrm{Ker}}}
\newcommand{\Ext}{\mathrm{Ext}}
\newcommand{\Hom}{\mathrm{Hom}}
\newcommand{\YD}{\mathcal{YD}}
\newcommand{\dotrtimes}{\mathop{\raisebox{0.2ex}{\makebox[0.86em][l]{${\scriptstyle >\joinrel\lessdot}$}}\raisebox{0.12ex}{$\shortmid$}}}
\newcommand{\blackrtimes}{\mathop{\raisebox{0.2ex}{${\scriptstyle >\joinrel\blacktriangleleft}$}}}
\newcommand{\ad}{\mathrm{ad}}
\newcommand{\diag}{\mathrm{diag}}
\newcommand{\ord}{\mathrm{ord}}
\newcommand{\rank}{\mathop{\mathrm{rank}}}

\begin{document}

\title[quantized enveloping algebras]{Abelian and non-abelian second cohomologies of quantized enveloping algebras}
\author{Akira Masuoka}  
\address{Institute of Mathematics, University of Tsukuba, Tsukuba, Ibaraki 305-8571 Japan}
\email{akira@math.tsukuba.ac.jp}

\begin{abstract}
For a class of pointed Hopf algebras including the quantized enveloping algebras, we discuss cleft extensions, cocycle deformations and the second cohomology.
We present such a non-standard method of computing the abelian second cohomology that derives information from the non-abelian second cohomology classifying cleft extensions.
As a sample computation, a quantum analogue of Whitehead's second lemma for Lie-algebra cohomology is proved.
\end{abstract}
\keywords{Hopf algebra, quantized enveloping algebra, cleft extension, the second cohomology, cocycle deformation}
\subjclass[2000]{16W30, 17B37}

\maketitle

\section*{Introduction} 

Let us work over a field $k$, for a moment.
Let $H$ be a Hopf algebra (over $k$).
By the abelian 2nd cohomology $H^{2}(H,M)$ of $H$, we mean the 2nd $\Ext$ group $\Ext_{H}^{2}(k,M)$, where $M$ is a left, say, $H$-module;
it is canonically isomorphic to the 2nd Hochschild cohomology $HH^{2}(H,M)$, in which $M$ is regarded as trivial as a right $H$-module.

On the other hand, by a non-abelian $2$-cocycle, we mean a {\em (right) $2$-cocycle} $\sigma : H\otimes H\rightarrow k$, which is by definition a convolution-invertible $k$-linear map satisfying the non-abelian, multiplicative $2$-cocycle conditions given by (1.1), (1.2) in the text.
Given such $\sigma$, the product on $H$ is deformed in the following two ways
\begin{eqnarray*}
a\cdot b & := & \sum\sigma(a_{1},b_{1})a_{2}b_{2}, \\
a\cdot b & := & \sum\sigma(a_{1},b_{1})a_{2}b_{2}\sigma^{-1}(a_{3},b_{3}),
\end{eqnarray*}
where $a,b\in H$.
Let ${}_{\sigma}H$, $H^{\sigma}$ denote the algebras defined by these deformed products, respectively.
Then, ${}_{\sigma}H$ together with the original regular coaction ${}_{\sigma}H\rightarrow {}_{\sigma}H\otimes H$ form a right $H$-cleft Galois extension over $k$ \cite{D1}, which we will call a {\em right $H$-cleft object}.
As was shown by Doi \cite{D1}, the natural equivalence classes of all $2$-cocycles are in 1-1 correspondence with the isomorphism classes, $\Cleft(H)$, of all right $H$-cleft objects;
in the text, the last symbol will denote the isomorphism classes of all left $H$-cleft objects, though the left and the right $H$-cleft objects are in 1-1 correspondence.
On the other hand, $H^{\sigma}$ together with the original coalgebra structure form a Hopf algebra, which is called a {\em cocycle deformation} of $H$ \cite{D2}.
${}_{\sigma}H$ is a left $H^{\sigma}$-cleft object with respect to the original coproduct ${}_{\sigma}H\rightarrow H^{\sigma}\otimes {}_{\sigma}H$, so that it is an $(H^{\sigma},H)$-bicleft object.
By Schauenburg's biGalois theory \cite{S1}, $H^{\sigma}$ is characterized as the Hopf algebra $L$ with which ${}_{\sigma}H$ is an $(L,H)$-biGalois object.
In \cite{G1}, \cite{G2}, \cite{S2}, \cite{Bi}, the authors classified right $H$-cleft objects for various $H$, and determined their associated Hopf algebras $L$; see also \cite{M1}.

Let $U_{q}$ denote the quantized enveloping algebra, due to Drinfeld and Jimbo, associated to a symmetrizable GCM (generalized Cartan matrix) $\mathbb{A}$, in which $q$ may be a root of $1$; see Section 7.
It has canonical generators $X_{i}^{\pm}$ ($i=1,2,\ldots,n$) of skew-primitives.
Let $(U_{q})^{0}$ denote the graded Hopf algebra defined in the same way as $U_{q}$, except that $X_{i}^{+}$ and $X_{i}^{-}$ commute each other in $(U_{q})^{0}$.
As part of the main results of this paper,
\begin{enumerate}
\renewcommand{\labelenumi}{(\arabic{enumi})}
\item (Theorem 7.8.2) we prove that $U_{q}$ and $(U_{q})^{0}$ are cocycle deformations of each other, giving an explicit isomorphism $U_{q}\simeq(U_{q})^{0}$ of coalgebras, and
\item (Theorem 7.14) we classify the central $U_{q}$-cleft extensions over a given commutative algebra $Z\neq 0$, assuming $\det\mathbb{A}\neq 0$ (e.g., $\mathbb{A}$ is of finite type), and that the order of $q$ is large enough.
\end{enumerate}
By a {\em central $H$-cleft extension} over $Z$, we mean an $H$-cleft Galois extension over $Z$ in which $Z$ is central.
It is identified with an $H\otimes Z$-cleft object over $Z$;
this fact would suggest us to work over a commutative ring rather than a field.
Notice from (1) that the central $U_{q}$-cleft extensions and the central $(U_{q})^{0}$-cleft extensions both over a fixed $Z$ are in a categorical 1-1 correspondence.
Assuming that the GCM $\mathbb{A}$ is of finite type, G\"unther \cite{G1} classified just as in (2) above, the central $U_{q}$-cleft extensions, but his result, Satz 4.9, contains some minor error; see Remark 7.14.
His method of computing is shown in \cite{G2} in the $\mathfrak{sl}_{2}$-case.
Under the same assumption as in \cite{G1}, Kassel and Schneider \cite{KS} proved the result in (1).
Our method is less technical, describing explicitly a $((U_{q})^{0},U_{q})$-bicleft object $A$, say.
We remark that $U_{q}$ and $A$ are naturally filtered so that in particular in $U_{q}$, the zero part $(U_{q})_{0}$ is spanned by the grouplikes.
We will prove that the associated graded algebra $\gr A$ is a $((U_{q})^{0},\gr U_{q})$-bicleft object, and indeed $\gr A=(U_{q})^{0}$.
As an advantage of our method, the next result follows.
\begin{enumerate}
\renewcommand{\labelenumi}{(3)}
\item (Theorem 7.8.1) $(U_{q})^{0}$ and $\gr U_{q}$ are canonically isomorphic.
\end{enumerate}

Let $H$ be a Hopf algebra again in general.
Let $M$ be a vector space, and regard it as a trivial $(H,H)$-bimodule.
Let
\[ k_{M}=k\oplus M=TM/(T^{2}M) \]
denote the tensor algebra $TM$ divided by all homogeneous components of degree $\geq 2$.
This has the projection $\pi : k_{M}\rightarrow k$ as an augmentation.
A central $H$-cleft extension $A$ over $k_{M}$ is said to be {\em augmented}, if it has an augmentation $A\rightarrow k$ extending $\pi$.
Let $\ZCleft_{\varepsilon}(H;k_{M})$ denote the set of the isomorphism (preserving the augmentation) classes of those augmented extensions.
Given a Hochschild $2$-cocycle $s : H\otimes H\rightarrow M$, then
\[ \sigma : H\otimes H\rightarrow k_{M},\quad \sigma(a,b)=\varepsilon(a)\varepsilon(b)+s(a,b) \]
is a $2$-cocycle with values in $k_{M}$, and ${}_{\sigma}(H\otimes k_{M})$ is a central $H$-cleft extension over $k_{M}$ which is augmented with respect to $\varepsilon\otimes\pi : {}_{\sigma}(H\otimes k_{M})\rightarrow k$.
We will prove (see Proposition 1.9) that $s\mapsto {}_{\sigma}(H\otimes k_{M})$ gives rise to a bijection
\[ H^{2}(H,M)\xrightarrow{\simeq}\ZCleft_{\varepsilon}(H;k_{M}). \]
(We remark again that we will work on the opposite side in the text.)
The dual result with non-trivial coefficients was used in \cite{M2}, \cite{MO} to prove the equivariant smoothness of hereditary Hopf algebras, in which case information on $\ZCleft_{\varepsilon}(H;k_{M})$ was derived from $H^{2}(H,M)$.
In this paper we apply the bijection above in the opposite direction.
When $H=U_{q}$, we will derive information from the result in (2) above, to prove
\begin{enumerate}
\renewcommand{\labelenumi}{(4)}
\item (Theorem 7.16) $H^{2}(U_{q},M)=0$ under the same assumption as in (2) and the additional assumption $\ch k\neq 2$.
\end{enumerate}
This is an essential part of a quantum analogue of Whitehead's 2nd lemma for Lie-algebra cohomology; see Remark 7.16.2.

Let $k\neq 0$ be a commutative ring.
Throughout we work over $k$, which will be supposed to a field only in the last Section 7 and Appendix.
Then central cleft extensions can be captured as cleft objects, as was seen above.
On the contrary we have the disadvantage that Schauenburg's biGalois theory \cite{S1}, in which the relevant Hopf algebras are mostly supposed to be flat over $k$, cannot apply, since those Hopf algebras which we will treat with do not necessarily satisfy the assumption.
But, this does not matter, since, as we will see in Section 1, the necessary results from \cite{S1} hold true without the flatness assumption, in our restricted situation that the involved Galois objects are all cleft.

We will actually derive the results cited above in (1)--(4), from the corresponding results for more general Hopf algebras, $H^{\lambda}$, $H^{0}$, which generalize $U_{q}$, $(U_{q})^{0}$, respectively.
Adopting the standpoint of Andruskiewitsch and Schneider \cite{AS1,AS2,AS3}, we construct $H^{\lambda}$, $H^{0}$ from a braided $k$-module.
Let $V=\bigoplus_{i\in I}kx_{i}$ be a free $k$-module with free basis $(x_{i})_{i\in I}$ indexed by a finite set $I\neq\emptyset$.
Let
\[ c : V\otimes V\xrightarrow{\simeq}V\otimes V,\quad c(x_{i}\otimes x_{j})=q_{ij}x_{j}\otimes x_{i} \]
denote the diagonal-type braiding on $V$ which is given by a matrix $(q_{ij})_{i,j\in I}$ with entries $q_{ij}$ in the group $k^{\times}$ of units in $k$.
Suppose that $(V,c)$ is realized in a natural way as an object in the braided monoidal category ${}^{\Gamma}_{\Gamma}\YD$ of Yetter-Drinfeld modules over the group Hopf algebra $k\Gamma$, where $\Gamma$ is some abelian group.
Then the tensor algebra $TV$ is naturally a graded, braided Hopf algebra in ${}^{\Gamma}_{\Gamma}\YD$, which gives rise to an ordinary graded Hopf algebra $TV\dotrtimes\Gamma$, by the Radford-Majid bosonization.
Given a system $\lambda=(\lambda_{ij})_{i,j}$ of parameters $\lambda_{ij}$ in $k$, where $(i,j)$ are in a certain subset of $I\times I$, we define $H^{\lambda}$ as a quotient (ungraded) Hopf algebra of $TV\dotrtimes\Gamma$; see Definition 4.1.
In the special case when all $\lambda_{ij}$ are zeros, $H^{\lambda}$ is denoted by $H^{0}$ (see Definition 2.2);
it is of the form $S\dotrtimes\Gamma$, where $S$ is a quotient graded braided Hopf algebra of $TV$.
In Section 2, a left $H^{0}$-comodule algebra $A(\lambda)$ is constructed, and is proved to be an $H^{0}$-cleft object.
Moreover, it is proved to be an $(H^{0},H^{\lambda})$-bicleft object, in Section 4.
As a consequence we have that $H^{0}$ and $H^{\lambda}$ are cocycle deformations of each other, and also that $H^{0}$ and $\gr H^{\lambda}$ are canonically isomorphic; see Theorem 4.3, Proposition 4.4.
The $H^{0}$-cleft object $A(\lambda)$ does not concern any group $2$-cocycle, or more precisely, the $2$-cocycle $H^{0}\otimes H^{0}\rightarrow k$ associated to $A(\lambda)$ restricts to the trivial group $2$-cocycle $\Gamma\times\Gamma\rightarrow k^{\times}$.
Deformations of $H^{0}$ and $A(\lambda)$ by group $2$-cocycles, or more precisely by those $2$-cocycles $H^{0}\otimes H^{0}\rightarrow k$ which arise from group $2$-cocycles $\Gamma\times\Gamma\rightarrow k^{\times}$ along the projection $H^{0}\rightarrow k\Gamma$, are discussed in Section 3.
In Section 5, augmented central $H^{\lambda}$-cleft extensions over $k_{M}$ are discussed.
We remark that their isomorphism classes (or the cohomology $H^{2}(H^{\lambda},M))$ depend on $\lambda$, while we have $\Cleft(H^{0})\simeq\Cleft(H^{\lambda})$, since $H^{0}$ and $H^{\lambda}$ are cocycle deformations of each other.
In Section 6, we put the additional assumptions (C1)--(C4), which includes the assumption that the matrix $(q_{ij})$ is {\em of Cartan type} \cite{AS1} in the sense that
\[ q_{ij}q_{ji}=q_{ii}^{a_{ij}}\ \ \mbox{for all}\ \ i,j\in I \]
for some GCM $(a_{ij})_{i,j\in I}$.
Under the assumptions, the set $\Cleft(H^{\lambda})$ and the cohomology $H^{2}(H^{\lambda},M)$ are determined by Theorems 6.3, 6.4.
In Section 7, the results obtained so far are applied to prove results for the quantized enveloping algebra $U_{q}$ and the Borel subalgebra $B_{q}$.
In the Appendix, which was added on the occasion of revision, we refine Grunenfelder and Mastnak's proof of their very recent result \cite{GM}, which essentially asserts that the finite-dimensional pointed Hopf algebra $u(\mathcal{D},\lambda,\mu)$ due to Andruskiewitsch and Schneider \cite{AS3} is a cocycle deformation of the associated graded Hopf algebra.

\noindent
{\it Notation}. The unadorned $\otimes$ denotes the tensor product over $k$.
For a commutative algebra $Z$, we let $Z^{\times}$ denote the group of all units in $Z$.
For sets $J$, $X$, we let $X^{J}$ denote the set of all maps $J\rightarrow X$.
An element in $X^{J}$ is expressed so as $\alpha=(\alpha_{j})_{j\in J}$, regarded as a system of parameters $\alpha_{j}$ in $X$.

\section{Preliminaries}

\subsection{Cocycle deformations}

Let $H$ be a Hopf algebra, whose coalgebra structure will be denoted, as usual, by
\[ \Delta : H\rightarrow H\otimes H,\ \ \Delta(a)=\sum a_{1}\otimes a_{2};\quad \varepsilon : H\rightarrow k. \]
By an $H$-comodule we mean a {\em left} $H$-comodule $V$, whose structure will be denoted by $v\mapsto \sum v_{-1}\otimes v_{0}$, $V\rightarrow H\otimes V$.
The $H$-comodules naturally form a $k$-linear monoidal category, which we denote by ${}^{H}\mathcal{M}$.
This is not necessarily abelian since we do not assume that $H$ is flat over $k$.

A {\em (right) $2$-cocycle} for $H$ is a convolution-invertible $k$-linear map $H\otimes H\rightarrow k$ such that
\begin{equation}
\sum\sigma(a_{1},b_{1})\sigma(a_{2}b_{2},c)=\sum\sigma(b_{1},c_{1})\sigma(a,b_{2}c_{2})
\end{equation}
\begin{equation}
\sigma(a,1)=\varepsilon(a)=\sigma(1,a)
\end{equation}
for all $a,b,c\in H$.
Given such $\sigma$, let $H^{\sigma}$ denote the coalgebra $H$ endowed with the deformed product defined by
\[ a\cdot b:=\sum\sigma(a_{1},b_{1})a_{2}b_{2}\sigma^{-1}(a_{3},b_{3}), \]
where $a,b\in H$.
This $H^{\sigma}$ forms a Hopf algebra with respect to the original unit $1$ and the deformed antipode $S^{\sigma}$ as given in \cite[Theorem 1.6 (b)]{D2}.
Given $V\in {}^{H}\mathcal{M}$, we can regard it as an $H^{\sigma}$-comodule in the obvious way since $H^{\sigma}=H$ as a coalgebra.
Denote this $H^{\sigma}$-comodule by ${}_{\sigma}V$.
Then we have a $k$-linear monoidal isomorphism
\begin{equation}
{}^{H}\mathcal{M}\xrightarrow{\simeq}{}^{H^{\sigma}}\!\!\mathcal{M},\ \ V\mapsto {}_{\sigma}V
\end{equation}
with respect to the monoidal structure
\begin{equation}
{}_{\sigma}V\otimes{}_{\sigma}W\xrightarrow{\simeq}{}_{\sigma}(V\otimes W),\ \ v\otimes w\mapsto \sum\sigma(v_{-1},w_{-1})v_{0}\otimes w_{0},
\end{equation}
\[ k\xrightarrow{\mathrm{id}}k={}_{\sigma}k\ \ \mbox{(identity)}. \]

We remark that the inverse $\sigma^{-1}$ is a $2$-cocycle for $H^{\sigma}$, and $(H^{\sigma})^{\sigma^{-1}}=H$.
An inverse of (1.3) is given by $V\mapsto {}_{\sigma^{-1}}V$.

\begin{definition}{\it \cite{D2}}
We say that two Hopf algebras $L,H$ are {\em cocycle deformations} of each other if $L\simeq H^{\sigma}$ for some $2$-cocycle $\sigma$ for $H$.
The condition is symmetric, as is seen from the remark above.
\end{definition}

\subsection{Cleft objects}

A (left) $H$-comodule algebra $A$ is said to be {\em cleft} \cite{DT}, if there exists a convolution-invertible $H$-colinear map $\phi : H\rightarrow A$.
This $\phi$ can be chosen so that $\phi(1)=1$.
In this case, $\phi$ is called a {\em section} \cite{D1}.
\begin{definition}
By an $H$-cleft object we mean a cleft $H$-comodule algebra $A$ such that $k\simeq {}^{\co H}A$, where 
\[ {}^{\co H}A=\{a\in A\;|\;\sum a_{-1}\otimes a_{0}=1\otimes a\}, \]
the subalgebra of $H$-coinvariants.
\end{definition}

A section $\phi : H\rightarrow A$ of an $H$-cleft object $A$ is necessarily bijective.
\begin{lemma}
An $H$-comodule algebra map $A\rightarrow A'$ between $H$-cleft objects is necessarily an isomorphism.
\end{lemma}
\begin{proof}
This is a special case of \cite[Lemma 1.3]{M1}.
\end{proof}

Let
\[ \Cleft(H) \]
denote the set of the isomorphism classes of all $H$-cleft objects.

Let $\sigma : H\otimes H\rightarrow k$ be a $2$-cocycle.
If $A$ is an $H$-comodule algebra (or an algebra-object in ${}^{H}\mathcal{M}$), then ${}_{\sigma}A$ is an $H^{\sigma}$-comodule algebra with respect to the deformed product
\[ a\cdot b := \sum\sigma(a_{1},b_{1})a_{2}b_{2}\quad (a,b\in A). \]

\begin{lemma}
The monoidal isomorphism (1.3) gives a bijection
\[ \Cleft(H)\simeq \Cleft(H^{\sigma}). \]
\end{lemma}
\begin{proof}
A section $\phi : H\rightarrow A$ gives a section $H^{\sigma}\rightarrow {}_{\sigma}A$.
\end{proof}

As a special right $H$-comodule algebra, a {\em right $H$-cleft object} is defined in the obvious way.
Let $L$ be another Hopf algebra.
An {\em $(L,H)$-bicleft object} is an $(L,H)$-bicomodule algebra which is a cleft object on both sides.
Let $A$ be such an object.
Then the cotensor product
\begin{equation}
A\Box_{H} : {}^{H}\mathcal{M}\rightarrow {}^{L}\mathcal{M}
\end{equation}
gives a $k$-linear monoidal functor whose monoidal structure is given by those natural isomorphisms
\begin{equation}
(A\Box_{H}V)\otimes(A\Box_{H}W)\xrightarrow{\simeq}A\Box_{H}(V\otimes W),
\end{equation}
\begin{equation}
k\xrightarrow{\simeq}A\Box_{H}k
\end{equation}
that arise from the algebra structure on $A$.
Working over a commutative ring, we should be careful about the cotensor product.
But, since $A\simeq H$ as right $H$-comodules, and since we have a natural isomorphism $V\simeq H\Box_{H}V$, it follows that the maps in (1.6), (1.7) are indeed isomorphisms in ${}^{L}\mathcal{M}$.

Let $\sigma : H\otimes H\rightarrow k$ be a $2$-cocycle.
Then, ${}_{\sigma}H$ is a right $H$-cleft object whose structure is given by the coproduct $\Delta : {}_{\sigma}H\rightarrow {}_{\sigma}H\otimes H$; see \cite[Proposition 1.1]{D1} or \cite[7.2.2]{Mo}.
Moreover, ${}_{\sigma}H$ is an $(H^{\sigma},H)$-bicleft object, by Lemma 1.4.
The monoidal functor ${}_{\sigma}H\Box_{H} : {}^{H}\mathcal{M}\rightarrow {}^{H^{\sigma}}\!\!\mathcal{M}$ as given by (1.5) is naturally isomorphic to the isomorphism (1.3).

Conversely, let $A$ be a right $H$-cleft object.
By \cite[Theorem 11]{DT} (or see \cite[7.2.2]{Mo}), a section $\phi : H\rightarrow A$ turns into an isomorphism of right $H$-cleft objects, if we replace $H$ with ${}_{\sigma}H$, where $\sigma : H\otimes H\rightarrow k$ is the $2$-cocycle defined by
\[ \sigma(a,b)=\sum\phi(a_{1})\phi(b_{1})\phi^{-1}(a_{2}b_{2})\quad (a,b\in H). \]

\begin{lemma}[{cf.\ \cite[Theorem 3.9]{S1}}]
Suppose in addition that $A$ is an $(L,H)$-bicleft object with $L$-coaction $\rho : A\rightarrow L\otimes A$.
Suppose that there is a coalgebra isomorphism $f : H\xrightarrow{\simeq}L$ which makes 
\[ \begin{CD}
H @>{\Delta}>> H\otimes H \\
@V{\phi}VV @VV{f\otimes\phi}V \\
A @>{\rho}>> L\otimes A
\end{CD} \]
commute, then $f$ turns into an isomorphism $H^{\sigma}\xrightarrow{\simeq}L$ of Hopf algebras, so that $L$ and $H$ are cocycle deformations of each other, and the monoidal functor (1.5) is an equivalence.
\end{lemma}
\begin{proof}
As is easily seen, in order for $\Delta : {}_{\sigma}H\rightarrow H\otimes {}_{\sigma}H$ to be an algebra map, the algebra structure on the coalgebra $H$ must coincide with that structure on $H^{\sigma}$.
This proves the lemma.
\end{proof}

A natural equivalence relation is defined among $2$-cocycles (see \cite{D1}), so that $\sigma \mapsto {}_{\sigma}H$ gives rise to a bijection from the set of the equivalence classes of $2$-cocycles for $H$ onto the set $\Cleft^{r}(H)$ of the isomorphism classes of right $H$-cleft objects.

A {\em left $2$-cocycle} $\tau : H\otimes H\rightarrow k$ is a convolution-invertible $k$-linear map which satisfies the mirror-image formulae of (1.1), (1.2).
Given such $\tau$, we have the left $H$-cleft object $H_{\tau}$ which is constructed just as ${}_{\sigma}H$, but the side is converted.
If $\sigma$ is a right $2$-cocycle, $\sigma^{-1}$ is a left $2$-cocycle, and vice versa.
One sees that ${}_{\sigma}H\mapsto H_{\sigma^{-1}}$ gives a bijection
\[ \Cleft^{r}(H)\simeq \Cleft(H). \]

\subsection{Augmented cleft extensions and the 2nd cohomology}

Let $Z\neq 0$ be a commutative algebra.
A {\em central $H$-cleft extension} over $Z$ is a cleft $H$-comodule algebra $A$ which includes $Z$ as a central subalgebra so that $Z={}^{\co H}A$.
Let
\begin{equation}
\ZCleft(H;Z)
\end{equation}
denote the set of the isomorphism classes of all central $H$-cleft extensions over $Z$.
It is naturally identified with the set $\Cleft(H\otimes Z)$ of the isomorphism classes of $H\otimes Z$-cleft objects over $Z$.

Suppose that $Z$ is augmented with respect to an algebra map $\varepsilon_{Z} : Z\rightarrow k$.

\begin{definition}
An {\em augmented central $H$-cleft extension} over $Z$ is a pair $(A,\varepsilon_{A})$ of a central $H$-cleft extension $A$ over $Z$, together with an algebra map $\varepsilon_{A} : A\rightarrow k$ extending $\varepsilon_{Z}$.
An {\em isomorphism} of those extensions is an augmented, $H$-comodule algebra isomorphism fixing each elements in $Z$.
We let
\begin{equation}
\ZCleft_{\varepsilon}(H;Z)
\end{equation}
denote the set of the isomorphism classes of augmented central $H$-cleft extensions over $Z$.
\end{definition}

Let $A=(A,\varepsilon_{A})$ be such an extension defined as above.
\begin{lemma}
A section $\phi : H\rightarrow A$ can be chosen so that it is augmented, or namely $\varepsilon = \varepsilon_{A}\circ\phi$.
\end{lemma}
\begin{proof}
Replace $\phi$ with $a\mapsto \sum\phi(a_{1})\varepsilon_{A}(\phi^{-1}(a_{2}))$, $H\rightarrow A$.
Notice that this is augmented, being still a section.
\end{proof}

A section $\phi : H\rightarrow A$ gives rise to a left $2$-cocycle with values in $Z$,
\[ \tau : H\otimes H\rightarrow Z,\quad \tau(a,b)=\sum\phi^{-1}(a_{1},b_{1})\phi(a_{2})\phi(b_{2}).  \]
Suppose that $\phi$ is augmented.
Then, $\tau$ is augmented in the sense $\varepsilon\otimes\varepsilon = \varepsilon_{Z}\circ\tau$, whence $\tau-\varepsilon\otimes\varepsilon$ has values in $\Ker\varepsilon_{Z}$.
We then have the central $H$-cleft extension $(H\otimes Z)_{\tau}$ over $Z$, which is augmented with respect to $\varepsilon\otimes\varepsilon_{Z} : (H\otimes Z)_{\tau}\rightarrow k$.
One sees that
\[ (H\otimes Z)_{\tau}\rightarrow A,\quad a\otimes z\mapsto \phi(a)z \]
is an isomorphism of augmented central $H$-cleft extensions over $Z$.

\begin{definition}{\it \cite[Chap.\ X]{CE}}
Given a right $H$-module $M$, we define the {\em $n$-th cohomology $k$-module} $H^{n}(H,M)$ of $H$ with coefficients in $M$ by
\[ H^{n}(H,M)=\Ext_{(H,k)}^{n}(k,M)\quad (n\geq 0), \]
where the right-hand side denotes the relative $n$-th $\Ext$ group \cite[Chap.\ IX]{Mac}, in which $k$ is regarded as a trivial right $H$-module through $\varepsilon$.
\end{definition}

Recall that $\Ext_{(H,k)}^{n}(k,M)$ is naturally isomorphic to the Hochschild cohomology $HH^{n}(H,{}_{\varepsilon}M)$, where ${}_{\varepsilon}M$ denote the $(H,H)$-bimodule $M$ with the original right $H$-action and the trivial left $H$-action through $\varepsilon$.
For the standard complexes for computing $HH^{n}(H,{}_{\varepsilon}M)$ gives a non-standard complex for computing $\Ext^{n}_{(H,k)}(k,M)$.
It follows that if $M$ is just a $k$-module, then it is possibly regarded as a trivial $H$-module both on left- and right-hand sides, but the $k$-module $\Ext_{(H,k)}^{n}(k,M)$ does not depend on the side.
Indeed, the (bi)modules over a Hopf algebra which we are going to treat with are all trivial, except in Remarks 1.10.2, 7.16.3.

Fix a $k$-module $M$, which we regard as a trivial $(H,H)$-bimodule.
Define an algebra $k_{M}$ by
\begin{equation}
k_{M}=k\oplus M=TM/(T^{2}M),
\end{equation}
the tensor algebra $TM$ divided by all homogeneous components of degree $\geq 2$.
Thus, $mm'=0$ in $k_{M}$, if $m,m'\in M$.
This is augmented with respect to the projection $k_{M}\rightarrow k$.
\begin{proposition}
(1) An augmented left $2$-cocycle $\tau : H\otimes H\rightarrow k_{M}$ arises uniquely from a (normalized) Hochschild $2$-cocycle $t : H\otimes H\rightarrow M$, so that
\begin{equation}
\tau(a,b)=\varepsilon(a)\varepsilon(b)+t(a,b)\quad (a,b\in H).
\end{equation}

(2) Given a Hochschild $2$-cocycle $t : H\otimes H\rightarrow M$, define a left $2$-cocycle $\tau : H\otimes H\rightarrow k_{M}$ by (1.11), which is augmented.
Then, $t\mapsto (H\otimes k_{M})_{\tau}$ gives rise to a bijection
\[ H^{2}(H,M)\xrightarrow{\simeq}\ZCleft_{\varepsilon}(H;k_{M}). \]
\end{proposition}
\begin{proof}
This last map is surjective, as is seen from the argument preceding Definition 1.8.
For the remaining, see the argument given in the proof of \cite[Theorem 4.1]{M2}.
\end{proof}

\begin{remark}
(1) As is seen from the proposition above, augmented left $2$-cocycles (with values in $k_{M}$) are precisely augmented right $2$-cocycles.
There is a bijection between the set $\ZCleft_{\varepsilon}(H;k_{M})$ and its right-hand side version, say $\ZCleft_{\varepsilon}^{r}(H;k_{M})$.

(2) {\em Augmented} (not necessarily central) {\em $H$-cleft extensions} over $k_{M}$ are defined in the obvious way.
Let $\Cleft_{\varepsilon}(H;k_{M})$ denote the set of their (augmented) isomorphism classes.
The last cited argument in \cite{M2} actually proves a bijection
\[ \bigsqcup_{\alpha}H^{2}(H,(M,\alpha))\simeq\Cleft_{\varepsilon}(H;k_{M}), \]
where $\alpha$ runs through all right $H$-module structures $\alpha : M\otimes H\rightarrow M$ on $M$.
See also \cite[Remark 1.6]{MO}.
\end{remark}

\subsection{Yetter-Drinfeld modules over a group}

Let $\Gamma$ be a group.
We do not assume here that $\Gamma$ is abelian, though we will do in later sections.

Let ${}^{\Gamma}_{\Gamma}\YD$ denote the $k$-linear braided monoidal category of Yetter-Drinfeld modules over the group Hopf algebra $k\Gamma$, which is abelian; see \cite[Chap.\ IX, Sect.\ 5]{K}, \cite[Sect.\ 10.6]{Mo}.
Thus an object in ${}^{\Gamma}_{\Gamma}\YD$ is a $\Gamma$-graded $k$-module $V=\bigoplus_{g\in\Gamma}V_{g}$ given a left $k\Gamma$-module structure such that
\[ gV_{h}\subset V_{ghg^{-1}}\quad (g,h\in\Gamma). \]
If $\Gamma$ is abelian, this condition means that each component $V_{g}$ is $\Gamma$-stable.
If an element $0\neq v\in V$ is homogeneous, and $v\in V_{g}$, then we write $|v|=g$.
Besides the obvious tensor product $V\otimes W$ and unit $k$, ${}^{\Gamma}_{\Gamma}\YD$ has the braiding
\[ c : V\otimes W\xrightarrow{\simeq}W\otimes V,\quad c(v\otimes w)=|v|w\otimes v. \]

By {\em cocycles, coboundaries} or {\em cochains} for $\Gamma$, we mean those in the normalized standard complex for computing the group cohomology $H^{\bullet}(\Gamma,\ \ )$.
The $\Gamma$-module of coefficients will be always trivial.
Fix a $2$-cocycle $\sigma : \Gamma\times\Gamma\rightarrow k^{\times}$.
It gives rise to a right (and at the same time, left) $2$-cocycle $k\Gamma\otimes k\Gamma\rightarrow k$ for the group Hopf algebra $k\Gamma$.
As a special case of (1.3), we have a $k$-linear monoidal isomorphism $V\mapsto {}_{\sigma}V$, ${}^{k\Gamma}\mathcal{M}\xrightarrow{\simeq}{}^{k\Gamma}\mathcal{M}$.
The structure (1.4) now turns into
\begin{equation}
{}_{\sigma}V\otimes {}_{\sigma}W\xrightarrow{\simeq} {}_{\sigma}(V\otimes W),\quad
v\otimes w\mapsto\sigma(|v|,|w|)v\otimes w.
\end{equation}
If $V\in{}^{\Gamma}_{\Gamma}\YD$, endow the $\Gamma$-graded $k$-module ${}_{\sigma}V$ with the deformed $\Gamma$-action
\begin{equation}
g\cdot v := \frac{\sigma(g,h)}{\sigma(ghg^{-1},g)}gv,
\end{equation}
where $g,h\in\Gamma$, $v\in V_{h}$; cf.\ \cite[Lemma 2.12]{AS2}.
One sees ${}_{\sigma}V\in {}^{\Gamma}_{\Gamma}\YD$.

\begin{proposition}[Takeuchi]
Moreover, $V\mapsto {}_{\sigma}V$ gives a $k$-linear monoidal isomorphism ${}^{\Gamma}_{\Gamma}\YD\xrightarrow{\simeq}{}^{\Gamma}_{\Gamma}\YD$ which preserves the braiding.
\end{proposition}

Let $R$ be a braided Hopf algebra in ${}^{\Gamma}_{\Gamma}\YD$.
Then, ${}_{\sigma}R$ is necessarily a braided Hopf algebra in ${}^{\Gamma}_{\Gamma}\YD$.

\begin{proposition}[Takeuchi]
Explicitly, ${}_{\sigma}R$ has the following deformed product and coproduct
\begin{eqnarray*}
a\cdot b & := & \sigma(g,h)ab\quad (a\in R_{g},\ b\in R_{h}), \\
\Delta(a) & := & \sum\sigma^{-1}(|a_{1}|,|a_{2}|)a_{1}\otimes a_{2}\quad (a\in R),
\end{eqnarray*}
while its unit, counit and antipode are the same as the original ones.
\end{proposition}

Keep $R$ as above.
By the Radford-Majid bosonization \cite{R} (or see \cite[Sect.\ 10.6]{Mo}), we have the ordinary Hopf algebra $R\dotrtimes\Gamma$.
As a $k$-module, $R\dotrtimes\Gamma = R\otimes k\Gamma$.
To denote an element, we will write $ag$ for $a\otimes g$, where $a\in R$, $g\in\Gamma$.
$R\dotrtimes\Gamma$ is the smash (or semi-direct) product $R\rtimes\Gamma$ as an algebra, and is the smash coproduct $R\blackrtimes k\Gamma$ as a coalgebra.
Explicitly, its product and coproduct are given by
\begin{eqnarray*}
(ag)(bh) & = & (a(gb))(gh), \\
\Delta(ag) & = & \sum a_{1}(g|a_{2}|)\otimes a_{2}g,
\end{eqnarray*}
where $a\in R$, $g,h\in\Gamma$.
From ${}_{\sigma}R$, we obtain the Hopf algebra ${}_{\sigma}R\dotrtimes\Gamma$.
Through the projection $R\dotrtimes\Gamma\rightarrow k\Gamma$, $\sigma$ gives rise to a right $2$-cocycle for $R\dotrtimes\Gamma$, which constructs the deformed $(R\dotrtimes\Gamma)^{\sigma}$.

\begin{proposition}[Takeuchi]
$ag\mapsto\sigma(|a|,g)ag$ gives a Hopf algebra isomorphism
\[ {}_{\sigma}R\dotrtimes\Gamma\xrightarrow{\simeq}(R\dotrtimes\Gamma)^{\sigma}. \]
\end{proposition}

The author learned the last three propositions, which are all directly verified, from M.\, Takeuchi's lectures in 2002 for the graduate course at Tsukuba.

\subsection{Graded coalgebras and Hopf algebras}

By {\em graded coalgebras} or {\em Hopf algebras}, we mean $\mathbb{N}$-graded those objects, where $\mathbb{N}=\{0,1,2,\ldots\}$.
Let $C=\bigoplus_{n=0}^{\infty}C(n)$ be a graded coalgebra.
Then, $C(0)$ is a coalgebra.
Being a direct summand of $C$, it can be called a subcoalgebra of $C$.

\begin{lemma}
Let $A$ be an (ungraded) algebra.
A $k$-linear map $f : C\rightarrow A$ is convolution-invertible if and only if the restriction $f|_{C(0)} : C(0)\rightarrow A$ onto $C(0)$ is.
\end{lemma}
\begin{proof}
To see the essential `if' part, notice that if $f|_{C(0)}$ is convolution-invertible, then for every $n$, the restriction $f|_{C_{n}}$ onto $C_{n}:=\bigoplus_{i=0}^{n}C(i)$ is, too, since the kernel of the restriction map $\Hom(C_{n},A)\rightarrow\Hom(C(0),A)$ is a nilpotent ideal.
\end{proof}

\begin{corollary}
Let $H$ be a bialgebra.
Suppose that $H$ is graded as a coalgebra so that the zero component $H(0)$ is a group Hopf algebra. Then, $H$ has a bijective antipode.
\end{corollary}
\begin{proof}
By Lemma 1.14, the identity maps $H\rightarrow H$, $H^{\mathrm{cop}}\rightarrow H$ are both convolution-invertible, since their restrictions onto the zero components are.
By \cite[Proposition 7]{DT}, $H$ has an antipode and a pode, which are necessarily mutual composite-inverses.
\end{proof}

Those Hopf algebras which we are going to treat with all satisfy the assumption of Corollary 1.15, whence their antipodes are bijective, by the corollary.

\section{The Hopf algebra $H^{0}$ and its cleft objects $A(\lambda)$}

Let us list up the symbols which will be defined in this section, and will be used throughout in what follows.
\begin{center}
$\Gamma$, $\hat{\Gamma}$, $I$ ($=I_{1}\sqcup I_{2}\sqcup\cdots\sqcup I_{N}$) \vspace{5pt} \\
$V=\{x_{i},g_{i},\chi_{i}\}_{i\in I}$, $q=(q_{ij})_{i,j\in I}$ \vspace{5pt} \\
$\Theta$, $\Xi$, $\mathcal{R}$, $\mathcal{R}^{0}$, $\mathcal{R}(\lambda)$, $H^{0}$, $A(\lambda)$
\end{center}

Let $\Gamma$ be an abelian group, and let $\hat{\Gamma}$ denote the character group, i.e., the group of all homomorphisms $\Gamma\rightarrow k^{\times}$.
Let $I$ be a non-empty finite set.

Let $q=(q_{ij})_{i,j\in I}$ be a square matrix with entries $q_{ij}$ in $k^{\times}$.
Let $g_{i}\in\Gamma$, $\chi_{i}\in\hat{\Gamma}$, such that
\[ \chi_{j}(g_{i})=q_{ij}\quad (i,j\in I). \]
Let $V=\bigoplus_{i\in I}kx_{i}$ be a free $k$-module with basis $(x_{i})_{i\in I}$.
This is made into an object in ${}^{\Gamma}_{\Gamma}\YD$ with respect to the structure
\[ |x_{i}|=g_{i},\ \ gx_{i}=\chi_{i}(g)x_{i}\quad (g\in\Gamma). \]
To indicate the structure we will write
\begin{equation}
V=\{x_{i},g_{i},\chi_{i}\}_{i\in I}\in {}^{\Gamma}_{\Gamma}\YD.
\end{equation}
The associated braiding is given by
\[ c : V\otimes V\xrightarrow{\simeq}V\otimes V,\quad c(x_{i}\otimes x_{j})=q_{ij}x_{j}\otimes x_{i}. \]
As a braided $k$-module, $V$ is thus {\em of diagonal type} \cite{AS1}, afforded by the matrix $q=(q_{ij})$.

The tensor algebra $TV$ on $V$ is naturally a graded, braided Hopf algebra in ${}^{\Gamma}_{\Gamma}\YD$, in which every element in $V$ is a primitive. Its bosonization will be denoted by
\begin{equation}
F=TV\dotrtimes\Gamma.
\end{equation}
This is a graded Hopf algebra with zero component $k\Gamma$.
Let $B_{l}$ denote the braid group of degree $l$.
The braiding $c$ on $V$ gives rise to a $kB_{l}$-module structure on the $l$-th power tensors $T^{l}V$.
Especially, the {\em braided commutator} \cite{AS1} is given by
\[ (\ad_{c}\,x_{i})(x_{j})=(1-c)(x_{i}\otimes x_{j})=x_{i}\otimes x_{j}-q_{ij}x_{j}\otimes x_{i}. \]

\begin{lemma}
If $q_{ij}q_{ji}=1$, then $(\ad_{c}\,x_{i})(x_{j})$ is a primitive in $TV$, and is a $(1,g_{i}g_{j})$-primitive in $F$.
\end{lemma}

This is well known.
In general an element $y$ in a coalgebra is called a {\em $(g,h)$-primitive}, where $g,h$ are grouplikes, if $\Delta(y)=y\otimes g+h\otimes y$; it is called a {\em skew-primitive} without specifying $g,h$.

Associated to $q=(q_{ij})$ is the graph having $I$ as the set of vertices, in which distinct vertices $i,j$ are connected by an edge if and only if $q_{ij}q_{ji}\neq 1$.
Fix a partition 
\[ I=I_{1}\sqcup I_{2}\sqcup\cdots\sqcup I_{N} \]
of $I$ into non-empty subsets $I_{r}$ ($r=1,2,\ldots,N$) which are mutually disconnected (in the graph) and disjoint.
We write $|i|=r$, if $i\in I_{r}$.
We have
\begin{equation}
q_{ij}q_{ji}=1\ \ \mbox{whenever}\ \ |i|\neq |j|.
\end{equation}
Define
\begin{equation}
\Theta = \{(i,j)\in I\times I\;|\;|i|>|j|\}.
\end{equation}

Let $\mathcal{R}$ be a subset of $TV$ consisting of elements of the form
\begin{equation}
\xi(x_{i_{1}}\otimes x_{i_{2}}\otimes\cdots\otimes x_{i_{l}})\quad (l\geq 2,\ \ \xi\in kB_{l}),
\end{equation}
where $|i_{1}|=|i_{2}|=\cdots=|i_{l}|$.
Notice that the element above is of degree $g_{i_{1}}g_{i_{2}}\cdots g_{i_{l}}$ in $\Gamma$.
Suppose that $\mathcal{R}$ includes a sequence of finite length $L\geq 0$, say, of subsets
\[ \emptyset=\mathcal{R}_{(0)}\subset\mathcal{R}_{(1)}\subset\dotsb\subset\mathcal{R}_{(L)}=\mathcal{R} \]
such that for each $0<h\leq L$, every element in $\mathcal{R}_{(h)}$ is a primitive modulo $(\mathcal{R}_{(h-1)})$ in $TV$; in particular, every element in $\mathcal{R}_{(1)}$ is required to be a primitive.
One then sees inductively that each ideal $(\mathcal{R}_{(h)})$ in $TV$ generated by $\mathcal{R}_{(h)}$ is a homogeneous braided Hopf ideal of $TV$.
Let $\mathcal{R}^{0}$ denote the union of the set $\mathcal{R}$ with the primitives (see Lemma 2.1)
\begin{equation}
(\ad_{c}\,x_{i})(x_{j})\quad ((i,j)\in\Theta).
\end{equation}
Obviously, $(\mathcal{R}^{0})$ is also a homogeneous braided Hopf ideal in $TV$.

\begin{definition}
Let
\[ S=TV/(\mathcal{R}^{0}) \]
denote the quotient graded braided Hopf algebra, and let
\[ H^{0}=S\dotrtimes\Gamma \]
denote the bosonization.
Notice that $H^{0}=F/(\mathcal{R}^{0})$, which denotes the quotient graded Hopf algebra of $F$ by the ideal $(\mathcal{R}^{0})$ generated by $\mathcal{R}^{0}$ in $F$.
\end{definition}

Define 
\begin{equation}
\Xi=\{(i,j)\in\Theta\;|\;\chi_{i}\chi_{j}=1\},
\end{equation}
where the $1$ in the condition denotes the constant map $g\mapsto 1$.

\begin{definition}
Given $\lambda=(\lambda_{ij})\in k^{\Xi}$, let $\mathcal{R}(\lambda)$ denote the union of the set $\mathcal{R}$ with the elements
\begin{equation}
(\ad_{c}\,x_{i})(x_{j})-\lambda_{ij}g_{i}g_{j}\quad ((i,j)\in\Xi),
\end{equation}
\begin{equation}
(\ad_{c}\,x_{i})(x_{j})\quad ((i,j)\in\Theta\setminus\Xi)
\end{equation}
in F. Define a quotient algebra of $F$ by 
\[ A(\lambda)=F/(\mathcal{R}(\lambda)). \]
\end{definition}

If $\lambda_{ij}=0$ for all $(i,j)\in\Xi$, then $\mathcal{R}(\lambda)=\mathcal{R}^{0}$, and $A(\lambda)=H^{0}$.
We unify the expressions (2.8), (2.9) into
\begin{equation}
(\ad_{c}\,x_{i})(x_{j})-\lambda_{ij}g_{i}g_{j}\quad ((i,j)\in\Theta),
\end{equation}
by setting $\lambda_{ij}=0$ for $(i,j)\in\Theta\setminus\Xi$.

\begin{lemma}
The coproduct $\Delta : F\rightarrow F\otimes F$ composed with the quotient map $F\otimes F\rightarrow H^{0}\otimes A(\lambda)$ factors through $A(\lambda)\rightarrow H^{0}\otimes A(\lambda)$, by which $A(\lambda)$ is an $H^{0}$-comodule algebra.
\end{lemma}
\begin{proof}
Notice first that the element $y$, say, given by (2.10) is sent by $\Delta$ to $(\ad_{c}\,x_{i})(x_{j})\otimes 1+g_{i}g_{j}\otimes y$. Then the result follows by induction on $L$.
\end{proof}

We have a direct sum decomposition
\[ V=V_{1}\oplus V_{2}\oplus\cdots\oplus V_{N}\ \ \mbox{in}\ \ {}^{\Gamma}_{\Gamma}\YD, \]
where $V_{r}=\bigoplus_{i\in I_{r}}kx_{i}$.
The set $\mathcal{R}$ is divided so as
\[ \mathcal{R}=\mathcal{R}_{1}\sqcup\mathcal{R}_{2}\sqcup\cdots\sqcup\mathcal{R}_{N} \]
into those mutually disjoint subsets $\mathcal{R}_{r}$ which consists of all elements (2.5) such that $r=|i_{1}|$ ($=|i_{2}|=\cdots =|i_{l}|$).
We have the graded algebras $TV_{r}/(R_{r})$.
Define a $k$-module by
\begin{equation}
E=\frac{TV_{1}}{(\mathcal{R}_{1})}\otimes\frac{TV_{2}}{(\mathcal{R}_{2})}\otimes\cdots\otimes\frac{TV_{N}}{(\mathcal{R}_{N})}\otimes k\Gamma.
\end{equation}
The tensor product of the canonical algebra maps $TV_{r}/(\mathcal{R}_{r})\rightarrow A(\lambda)$ ($r=1,2,\ldots,N$), composed with the product $A(\lambda)^{\otimes N}\rightarrow A(\lambda)$, induces a right $k\Gamma$-linear map,
\begin{equation}
e(\lambda) : E\rightarrow A(\lambda).
\end{equation}
Especially we have
\[ e:=e(0) : E\rightarrow H^{0}. \]
We remark that each $TV_{r}/(\mathcal{R}_{r})$ is a graded, braided Hopf algebra in ${}^{\Gamma}_{\Gamma}\YD$.
As the bosonization of their braided tensor product, $E$ has a natural graded Hopf algebra structure.
We see easily that $e$ is a graded Hopf algebra map.

\begin{proposition}
(1) These maps $e(\lambda)$, $e$ are bijective.

(2) The composite
\[ \phi(\lambda) := e(\lambda)\circ e^{-1} : H^{0}\rightarrow A(\lambda) \]
is a section.

(3) $A(\lambda)$ is an $H^{0}$-cleft object.
\end{proposition}

\begin{remark}
Fix $1\leq r\leq N$.
By Lemma 2.4 and Part 1 above, $e(\lambda)$ (resp., $e$) embeds the Hopf algebra
\begin{equation}
E_{r}:=\frac{TV_{r}}{(\mathcal{R}_{r})}\dotrtimes\Gamma
\end{equation}
obtained as the bosonization of $TV_{r}/(\mathcal{R}_{r})$, into $A(\lambda)$ as a subalgebra (resp., into $H^{0}$ as a Hopf subalgebra), so that the structure map $A(\lambda)\rightarrow H^{0}\otimes A(\lambda)$ restricts to the coproduct of $E_{r}$.
Moreover, the section $\phi(\lambda)$ above restricts to the identity map of $E_{r}$.
In this sense, the pair $(A(\lambda),\phi(\lambda))$ restricts to the $E_{r}$-cleft object $E_{r}$ with the obvious section $\mathrm{id}$.
\end{remark}

We will prove the proposition by induction on $L$.

\begin{lemma}
Proposition 2.5 holds true in $L=0$.
\end{lemma}

For the proof of Part 1, we will apply Bergman's diamond lemma \cite[Proposition 7.1]{Be} for $R$-rings.
Recall that for an algebra $R$, an {\em $R$-ring} is an algebra given an algebra map from $R$.
\begin{proof}[Proof of Lemma 2.7]
Suppose $L=0$, or namely $\mathcal{R}=\emptyset$.
Then, $E=TV_{1}\otimes\dotsb\otimes TV_{N}\otimes k\Gamma$, in particular.

(1) It suffices to prove only for $e(\lambda)$.
Notice that $A(\lambda)$ is a $k\Gamma$-ring, which is generated by the elements $x_{i}$ ($i\in I$), and is defined by the relations
\begin{eqnarray}
x_{i}x_{j} & = & q_{ij}x_{j}x_{i}+\lambda_{ij}g_{i}g_{j}\quad ((i,j)\in\Theta), \\
gx_{i} & = & \chi_{i}(g)x_{i}g\quad (i\in I,\ g\in\Gamma),
\end{eqnarray}
where by convention, $\lambda_{ij}=0$ if $(i,j)\not\in\Xi$.
To apply \cite[Proposition 7.1]{Be}, we introduce a total ordering into $(x_{i})_{i\in I}$ (or into $I$), such that $i>j$ if $|i|>|j|$, and consider the lexicographic ordering on the monomials in $x_{i}$.
Then the relations (2.14), (2.15) are naturally regarded as a reduction system.
Essential is to prove that the overlap ambiguities which arise when we reduce
\begin{equation}
x_{i}x_{j}x_{k}\quad (|i|>|j|>|k|)
\end{equation}
are resolvable.
One sees by computation that the reduction results of (2.16) in two ways differ by
\[ (q_{ij}q_{ik}\lambda_{jk}-\lambda_{jk})x_{i}+(q_{ij}\lambda_{ik}-q_{jk}\lambda_{ik})x_{j}+(\lambda_{ij}-q_{ik}q_{jk}\lambda_{ij})x_{k}. \]
But, this is zero, implying the desired resolvability.
For if $i\in I$, $(i,k)\in\Theta$,
\begin{eqnarray*}
q_{ij}q_{ik}\lambda_{jk} & = & \chi_{i}\chi_{k}(g_{i})\lambda_{jk}=\left\{\begin{array}{ll} \lambda_{jk} & \mbox{if $(j,k)\in\Xi$} \\ 0 & \mbox{otherwise} \end{array}\right. \\
 & = & \lambda_{jk}.
\end{eqnarray*}
By \cite[Proposition 7.1]{Be}, $A(\lambda)$ has a right $k\Gamma$-free basis $x_{i_{1}}x_{i_{2}}\cdots x_{i_{l}}$, where $l\geq 0$, $|i_{1}|\leq |i_{2}|\leq\cdots\leq |i_{l}|$.
This proves Part 1.

(2) $\phi(\lambda)$ preserves unit, since $e(\lambda)$, $e$ do.
We have the obvious embedding $E\hookrightarrow F$ of (right $k\Gamma$-module) coalgebras, and its composite with $F\rightarrow A(\lambda)$ is precisely $e(\lambda)$.
We see from Lemma 2.4 that $\phi(\lambda)$ is $H^{0}$-colinear.
It is convolution-invertible by Lemma 1.14.

(3) Since $\phi(\lambda)$ is an $H^{0}$-colinear isomorphism, $k$ is precisely the $H^{0}$-coinvariants in $A(\lambda)$, whence Part 3 follows.
\end{proof}

Suppose that we have proved the proposition in $L-1$, and let $\tilde{H}^{0}$, $\tilde{A}(\lambda)$ denote those $H^{0}$, $A(\lambda)$ defined by supposing $\mathcal{R}=\mathcal{R}_{(L-1)}$.
They both then include
\begin{equation}
\tilde{T}_{r}:=TV_{r}/(\mathcal{R}_{r}\cap\mathcal{R}_{(L-1)})\quad (1\leq r\leq N),
\end{equation}
as well as $k\Gamma$, as subalgebras.
Choose an arbitrary element in $\mathcal{R}$ as given by (2.5), and let $y$ denote its image in $\tilde{T}_{r}$, where we have set $r=|i_{1}|$ ($=|i_{2}|=\dotsb=|i_{l}|$).
Notice that $\tilde{T}_{r}$ is graded so as $\tilde{T}_{r}=\bigoplus_{m=0}^{\infty}\tilde{T}_{r}(m)$, where $\tilde{T}_{r}(m)$ denotes the component of degree $m$.
Then, $y\in\tilde{T}_{r}(l)$.

\begin{lemma}
Keep the assumptions and the notation as above.
Suppose $|j|\neq r$.
Then both in $\tilde{H}^{0}$, and in $\tilde{A}(\lambda)$,
\[ x_{j}y=q_{ji_{1}}q_{ji_{2}}\cdots q_{ji_{l}}yx_{j}. \]
\end{lemma}
\begin{proof}
In $\tilde{H}^{0}$, this holds true since $y$ is a $k$-linear combination of the monomial $x_{i_{1}}x_{i_{2}}\cdots x_{i_{l}}$ and its permutations.

To prove in $\tilde{A}(\lambda)$, let $\tilde{\rho} : \tilde{A}(\lambda)\rightarrow \tilde{H}^{0}\otimes\tilde{A}(\lambda)$ denote the structure.
Set 
\[ z=x_{j}y-q_{ji_{1}}\cdots q_{ji_{l}}yx_{j}. \]
By assumption, $y$ is a primitive in $\tilde{T}_{r}$.
By the induction hypothesis, we see from Remark 2.6 that $\tilde{\rho}(y)=y\otimes 1+g\otimes y$, where $g=g_{i_{1}}\cdots g_{i_{l}}$.
Since $x_{j}g=q_{ji_{1}}\cdots q_{ji_{l}}gx_{j}$, we see
\[ \tilde{\rho}(z)=z\otimes 1+gg_{j}\otimes z=gg_{j}\otimes z. \]
Again by the induction hypothesis above, $z=cgg_{j}$ for some $c\in k$.
But, $c$ must be zero, whence $z=0$, as desired, since we see by induction on $m$ that in $\tilde{A}(\lambda)$,
\begin{equation}
x_{j}\tilde{T}_{r}(m)\subset \tilde{T}_{r}(m)x_{j}+\tilde{T}_{r}(m-1)\Gamma\quad (m>0),
\end{equation}
provided $|j|\neq r$.
\end{proof}

\begin{proof}[Proof of Proposition 2.5]
To prove in $L$, keep the assumptions and the notation as above.
Fix $r$, and let $\mathfrak{a}_{r}$ denote the ideal of $\tilde{T}_{r}$ generated by the natural image of $\mathcal{R}_{r}$.
Notice from (2.18) that $\tilde{T}_{r}\tilde{T}_{s}\Gamma=\tilde{T}_{s}\tilde{T}_{r}\Gamma$ in $\tilde{A}(\lambda)$, for all $s$.
By Lemma 2.8, we see that both in $\tilde{H}^{0}$ and in $\tilde{A}(\lambda)$, $\tilde{T}_{s}\mathfrak{a}_{r}\Gamma=\mathfrak{a}_{r}\tilde{T}_{s}\Gamma$ for all $s$.
Obviously, $\Gamma\mathfrak{a}_{r}=\mathfrak{a}_{r}\Gamma$.
Hence in the both algebras, the ideal generated by the image of $\mathcal{R}_{r}$ equals
\[ \tilde{T}_{1}\dotsm\tilde{T}_{r-1}\mathfrak{a}_{r}\tilde{T}_{r+1}\dotsm\tilde{T}_{N}\Gamma. \]
This proves Part 1.
We see that $\phi(\lambda)$ is $H^{0}$-colinear since it is induced from the $H^{0}$-colinear $\tilde{\phi}(\lambda)$.
The remaining can be proved as Lemma 2.7.
\end{proof}

\begin{remark}
Proposition 2.5 holds true even if we rearrange the $TV_{r}$ ($r=1,\ldots,N$) in $E$ in an arbitrary order.
To see this, it suffices in the proof of Lemma 2.7.1, to replace the relation (or reduction) (2.14) with the equivalent
\[ x_{j}x_{i}=q_{ji}x_{i}x_{j}-q_{ji}\lambda_{ij}g_{i}g_{j}, \]
if $|i|<|j|$ in the new ordering.
\end{remark}

\section{Deformations of $H^{0}$, $A(\lambda)$ by group $2$-cocycles}

Let $Z^{2}(\Gamma,k^{\times})$ denote the group of all normalized $2$-cocycles $\Gamma\times\Gamma\rightarrow k^{\times}$.
Let $\sigma\in Z^{2}(\Gamma,k^{\times})$.
We will write $\check{\sigma}$ for the inverse $\sigma^{-1}$.
From $V=\{x_{i},g_{i},\chi_{i}\}_{i\in I}$ in (2.1), we obtain by Proposition 1.11
\begin{equation}
{}_{\sigma}V=\{x_{i},g_{i},\chi_{i}^{\sigma}\}_{i\in I}\in {}^{\Gamma}_{\Gamma}\YD,
\end{equation}
where
\begin{equation}
\chi_{i}^{\sigma}(g)=\frac{\sigma(g,g_{i})}{\sigma(g_{i},g)}\chi_{i}(g)\quad (g\in\Gamma).
\end{equation}
The associated braiding $c^{\sigma} : {}_{\sigma}V\otimes {}_{\sigma}V\xrightarrow{\simeq} {}_{\sigma}V\otimes {}_{\sigma}V$ is given by $c^{\sigma}(x_{i}\otimes x_{j})=q_{ij}^{\sigma}x_{j}\otimes x_{i}$, where 
\begin{equation}
q_{ij}^{\sigma}=\frac{\sigma(g_{i},g_{j})}{\sigma(g_{j},g_{i})}q_{ij}.
\end{equation}
See \cite[Lemma 2.12]{AS2}.
Regarding $\sigma$ as a right $2$-cocycle for $H^{0}$ through the projection $H^{0}=S\dotrtimes\Gamma\rightarrow k\Gamma$, we obtain $(H^{0})^{\sigma}$, ${}_{\sigma}A(\lambda)$.
We aim to describe explicitly these deformed objects.

Choose another $\tau\in Z^{2}(\Gamma,k^{\times})$.
Then we have the $\Gamma$-graded algebra $TV\rtimes_{\tau}\Gamma$ of crossed product;
it differs from the smash product $TV\rtimes\Gamma$ only by
\begin{equation}
\bar{g}\bar{g}'=\tau(g,g')\overline{gg'}\quad (g,g'\in\Gamma),
\end{equation}
where $g$ ($=1\otimes g$) in $TV\rtimes\Gamma$ is denoted by $\bar{g}$ in $TV\rtimes_{\tau}\Gamma$.
Let $\mathcal{Q}$ be a set consisting of such elements in $TV\rtimes_{\tau}\Gamma$ of the form
\begin{equation}
\xi(x_{i_{1}}\otimes x_{i_{2}}\otimes\cdots\otimes x_{i_{l}})-c\bar{g}_{i_{1}}\bar{g}_{i_{2}}\cdots\bar{g}_{i_{l}},
\end{equation}
where $l\geq 2$, $\xi\in kB_{l}$, $c\in k$.
Notice that the element above is of degree $g_{i_{1}}g_{i_{2}}\cdots g_{i_{l}}$ in $\Gamma$.

\begin{definition}
Read every element (3.5) in $\mathcal{Q}$ as an element in $T{}_{\sigma}V\rtimes_{\sigma\tau}\Gamma$, by replacing the action by $\xi$ with what arises from $c^{\sigma}$, and taking $\bar{g}_{i_{1}}\bar{g}_{i_{2}}\cdots\bar{g}_{i_{l}}$ as the product given by $\sigma\tau$.
Let ${}_{\sigma}\mathcal{Q}$ denote the subset of $T{}_{\sigma}V\rtimes_{\sigma\tau}\Gamma$ consisting of thus read elements.
\end{definition}

\begin{lemma}
(1) The natural embeddings of ${}_{\sigma}V$, $\Gamma$ into ${}_{\sigma}(TV\rtimes_{\tau}\Gamma)$ induce an isomorphism of $\Gamma$-graded algebras,
\[ T{}_{\sigma}V\rtimes_{\sigma\tau}\Gamma\xrightarrow{\simeq} {}_{\sigma}(TV\rtimes_{\tau}\Gamma). \]

(2) The isomorphism above induces an isomorphism,
\[ T{}_{\sigma}V\rtimes_{\sigma\tau}\Gamma/({}_{\sigma}\mathcal{Q})\xrightarrow{\simeq} {}_{\sigma}(TV\rtimes_{\tau}\Gamma/(\mathcal{Q})). \]
\end{lemma}
\begin{proof}
(1) We see from (3.2) that such a $\Gamma$-graded algebra map as above is defined.
To see that it is an isomorphism, it suffices to apply the result below in the special case when $\mathcal{Q}=\emptyset$.

(2) In general let $A$ be a $\Gamma$-graded algebra including $V$ as a $\Gamma$-graded submodule.
Let
\begin{equation}
y = \xi(x_{i_{1}}\otimes x_{i_{2}}\otimes \cdots \otimes x_{i_{l}})\in T^{l}V,
\end{equation}
where $l\geq 2$, $\xi\in kB_{l}$, and let ${}_{\sigma}y$ be the corresponding element in $T^{l}{}_{\sigma}V$ read as above.
By considering the $l$-fold iteration ${}_{\sigma}V\otimes\cdots\otimes {}_{\sigma}V\xrightarrow{\simeq} {}_{\sigma}(V\otimes\cdots\otimes V)$ of (1.12), we see that the product of ${}_{\sigma}y$ taken in ${}_{\sigma}A$ differs from the product of $y$ taken in $A$, only by scalar multiplication by
\begin{equation}
u := \sigma(g_{i_{1}},g_{i_{2}})\sigma(g_{i_{1}}g_{i_{2}},g_{i_{3}})\cdots\sigma(g_{i_{1}}\cdots g_{i_{l-1}},g_{i_{l}}).
\end{equation}
Suppose $A=TV\rtimes_{\tau}\Gamma/(\mathcal{Q})$.
Then it follows that the composite
\[ T{}_{\sigma}V\rtimes_{\sigma\tau}\Gamma\rightarrow {}_{\sigma}(TV\rtimes_{\tau}\Gamma)\rightarrow {}_{\sigma}A \]
induces
\[ f : B:=T{}_{\sigma}V\rtimes_{\sigma\tau}\Gamma/({}_{\sigma}\mathcal{Q})\rightarrow {}_{\sigma}A. \]
This is surjective, since ${}_{\check{\sigma}}f$, whose image contains the natural generators, is surjective; recall here $\check{\sigma}=\sigma^{-1}$.
By applying the result to $B$, we have an epimorphism
\[ h : A=TV\rtimes_{\tau}\Gamma/(\mathcal{Q})\rightarrow {}_{\check{\sigma}}B. \]
Since the composite ${}_{\sigma}h\circ f$ is the identity map of $B$, $f$ is a isomorphism.
\end{proof}

\begin{remark}
Let $y$, $u$ be as in (3.6), (3.7).
As is seen from the proof above, the canonical, braided Hopf algebra map $T{}_{\sigma}V\rightarrow {}_{\sigma}(TV)$ is an isomorphism, under which ${}_{\sigma}y\mapsto uy$.
\end{remark}

Recall from Section 2 the constructions
\[ S=TV/(\mathcal{R}^{0}),\ \ H^{0}=TV\dotrtimes\Gamma/(\mathcal{R}^{0}),\ \ A(\lambda)=TV\rtimes\Gamma/(\mathcal{R}(\lambda)). \]
Notice from Remark 3.3 that the set ${}_{\sigma}\mathcal{R}^{0}$ together with
\[ \emptyset={}_{\sigma}\mathcal{R}_{(0)}\subset{}_{\sigma}\mathcal{R}_{(1)}\subset\dotsb\subset{}_{\sigma}\mathcal{R}_{(L)}={}_{\sigma}\mathcal{R} \]
satisfy the same conditions as what was required to $\mathcal{R}$, so that $T_{\sigma}V/({}_{\sigma}\mathcal{R}_{(h)})\simeq {}_{\sigma}(TV/(\mathcal{R}_{(h)}))$ ($0\leq h\leq L$) as graded braided Hopf algebras. Moreover,
\[ T{}_{\sigma}V/({}_{\sigma}\mathcal{R}^{0}),\ \ T{}_{\sigma}V\dotrtimes\Gamma/({}_{\sigma}\mathcal{R}^{0}) \]
are a graded, braided Hopf algebra in ${}^{\Gamma}_{\Gamma}\YD$, and an ordinary graded Hopf algebra, respectively.
${}_{\sigma}\mathcal{R}(\lambda)$ is the union of the set ${}_{\sigma}\mathcal{R}$ with the elements
\[ x_{i}\otimes x_{j}-q_{ij}^{\sigma}x_{j}\otimes x_{i}-\lambda_{ij}\bar{g}_{i}\bar{g}_{j}\quad ((i,j)\in\Theta) \]
in $T{}_{\sigma}V\rtimes_{\sigma}\Gamma$.

Notice that $TV\rtimes_{\sigma}\Gamma$ is a $TV\dotrtimes\Gamma$-comodule algebra with respect to the algebra map
\begin{equation}
TV\rtimes_{\sigma}\Gamma\rightarrow (TV\dotrtimes\Gamma)\otimes(TV\rtimes_{\sigma}\Gamma)
\end{equation}
given by $x_{i}\mapsto x_{i}\otimes 1+g_{i}\otimes x_{i}$ ($i\in I$), $\bar{g}\mapsto g\otimes\bar{g}$ ($g\in\Gamma$).
This is indeed a cleft object for which the identity map
\begin{equation}
\mathrm{id} : TV\dotrtimes\Gamma\rightarrow TV\rtimes_{\sigma}\Gamma,
\end{equation}
with both identified with $TV\otimes k\Gamma$, is an obvious section.
Analogously, $T{}_{\sigma}V\rtimes_{\sigma}\Gamma$ is a $T{}_{\sigma}V\dotrtimes\Gamma$-cleft object.
This induces a $T{}_{\sigma}V\dotrtimes\Gamma/({}_{\sigma}\mathcal{R}^{0})$-comodule algebra,
\[ T{}_{\sigma}V\rtimes_{\sigma}\Gamma/({}_{\sigma}\mathcal{R}(\lambda)), \]
as is easily seen.
This is a cleft object, too, as will be seen from the next proposition.

\begin{proposition}
(1) $x_{i}\mapsto x_{i}$ ($i\in I$), $g\mapsto g$ ($g\in\Gamma$) give rise to isomorphisms
\begin{eqnarray}
T{}_{\sigma}V/({}_{\sigma}\mathcal{R}^{0}) & \xrightarrow{\simeq} & {}_{\sigma}S, \nonumber \\
T{}_{\sigma}V\dotrtimes\Gamma/({}_{\sigma}\mathcal{R}^{0}) & \xrightarrow{\simeq} & {}_{\sigma}S\dotrtimes\Gamma=(H^{0})^{\sigma}
\end{eqnarray}
of graded, braided Hopf algebras in ${}^{\Gamma}_{\Gamma}\YD$, and of graded Hopf algebras, respectively.

(2) $x_{i}\mapsto x_{i}$ ($i\in I$), $\bar{g}\mapsto g$ ($g\in\Gamma$) give rise to an algebra isomorphism
\[ T{}_{\sigma}V\rtimes_{\sigma}\Gamma/({}_{\sigma}\mathcal{R}(\lambda))\xrightarrow{\simeq} {}_{\sigma}A(\lambda), \]
which is compatible with (3.10) through the left coactions.
\end{proposition}
\begin{proof}
This follows by Lemma 3.2 and Proposition 1.13.
\end{proof}

Let $\sigma\in Z^{2}(\Gamma,k^{\times})$.
Define $s_{i}^{\sigma}\in\hat{\Gamma}$ by
\begin{equation}
s_{i}^{\sigma}(g)=\frac{\sigma(g,g_{i})}{\sigma(g_{i},g)}\quad (g\in\Gamma).
\end{equation}
One sees $s_{i}^{\sigma}=\chi_{i}^{\sigma}/\chi_{i}$, whence $s_{i}^{\sigma}\in\hat{\Gamma}$; see (3.2).
Define 
\begin{equation}
\Xi(\sigma)=\{(i,j)\in\Theta\;|\;\chi_{i}\chi_{j}=s_{i}^{\sigma}s_{j}^{\sigma}\}.
\end{equation}
The condition given here is equivalent to $\chi_{i}^{\check{\sigma}}\chi_{j}^{\check{\sigma}}=1$.
If $\sigma=1$, then $\Xi(\sigma)=\Xi$; see (2.7).
If $\sigma$ and $\sigma'$ in $Z^{2}(\Gamma,k^{\times})$ are cohomologous, then $s_{i}^{\sigma}=s_{i}^{\sigma'}$, whence $\Xi(\sigma)=\Xi(\sigma')$.

\begin{definition}
Given $\mu=(\mu_{ij})\in k^{\Xi(\sigma)}$, let $\mathcal{R}(\sigma,\mu)$ denote the union of the set $\mathcal{R}$ with the elements
\begin{equation}
x_{i}\otimes x_{j}-q_{ij}x_{j}\otimes x_{i}-\mu_{ij}\bar{g}_{i}\bar{g}_{j}\quad ((i,j)\in\Theta)
\end{equation}
in $TV\rtimes {}_{\sigma}\Gamma$, where we set $\mu_{ij}=0$ if $(i,j)\not\in\Xi(\sigma)$.
Define
\begin{equation}
A(\sigma,\mu)=TV\rtimes_{\sigma}\Gamma/(\mathcal{R}(\sigma,\mu)).
\end{equation}
\end{definition}

The algebra $A(\lambda)$ equals $A(1,\lambda)$ with the notation above.

\begin{proposition}
The algebra map (3.8) induces an algebra map
\[ A(\sigma,\mu)\rightarrow H^{0}\otimes A(\sigma,\mu), \]
by which $A(\sigma,\mu)$ is an $H^{0}$-cleft object.
\end{proposition}
\begin{proof}
Reconsider the construction of $A(\lambda)$, starting from ${}_{\check{\sigma}}V$, in which case the set $\Xi$ in (2.7) is $\Xi(\sigma)$.
Set $\mathcal{Q}^{0}={}_{\check{\sigma}}\mathcal{R}^{0}$.
Given $\mu\in\Xi(\sigma)$, let $\mathcal{Q}(\mu)$ denote the union of the set ${}_{\check{\sigma}}\mathcal{R}$ with the elements
\[ x_{i}\otimes x_{j}-q_{ij}^{\check{\sigma}}x_{j}\otimes x_{i}-\mu_{ij}g_{i}g_{j}\quad ((i,j)\in\Theta) \]
in $T{}_{\check{\sigma}}V\dotrtimes\Gamma$, where by convention, $\mu_{ij}=0$ if $(i,j)\not\in\Xi(\sigma)$.
Then we have the $T{}_{\check{\sigma}}V\dotrtimes\Gamma/(\mathcal{Q}^{0})$-cleft object $T{}_{\check{\sigma}}V\rtimes\Gamma/(\mathcal{Q}(\mu))$.
Apply ${}_{\sigma}(\ \ )$, and use the isomorphisms in Proposition 3.4.
Then we see that $A(\sigma,\mu)$ is an $H^{0}$-cleft object, since ${}_{\sigma}\mathcal{Q}^{0}=\mathcal{R}^{0}$, ${}_{\sigma}\mathcal{Q}(\mu)=\mathcal{R}(\sigma,\mu)$.
The section
\[ T{}_{\check{\sigma}}V\dotrtimes\Gamma/(\mathcal{Q}^{0})\rightarrow T{}_{\check{\sigma}}V\rtimes\Gamma/(\mathcal{Q}(\mu)) \]
given by  Proposition 2.5.2, applied with ${}_{\sigma}(\ \ )$, gives rise to a section,
\begin{equation}
\phi(\sigma,\mu) : H^{0}\rightarrow A(\sigma,\mu).
\end{equation}
\end{proof}

\begin{remark}
Fix $1\leq r\leq N$.
Recall from Remark 2.6 the embedding $E_{r}\hookrightarrow A(\lambda)$, but consider now this to be given for ${}_{\sigma}V,\mathcal{Q},\mu$ instead of $V,\mathcal{R},\lambda$.
With ${}_{\sigma}(\ \ )$ applied, it gives rise to an embedding
\begin{equation}
E_{r}(\sigma):=\frac{TV_{r}}{(\mathcal{R}_{r})}\rtimes_{\sigma}\Gamma\hookrightarrow A(\sigma,\mu),
\end{equation}
which coincides with the canonical algebra map.
It follows that $A(\sigma,\mu)$ together with the section $\phi(\sigma,\mu)$ above restrict to the $E_{r}$-cleft object $E_{r}(\sigma)$ with the obvious section (see (3.9)), in the same sense as in Remark 2.6.
\end{remark}

We wish to classify all $A(\sigma,\mu)$ into isomorphism classes.
Since $V\oplus k\Gamma$ in $TV\dotrtimes\Gamma$ is isomorphically mapped into $H^{0}$, we denote the images of the natural basis elements $x_{i}$ ($i\in I$), $g$ ($\in\Gamma$) in $V\oplus k\Gamma$ by the same symbols $x_{i}$, $g$.
Let $A$ be an $H^{0}$-cleft object.

\begin{definition}
An element $a$ in $A$ is called an {\em NB} ({\em normal basis}) {\em element} for $x_{i}$ (resp., $g$) if there is a section $\phi : H^{0}\rightarrow A$ such that $a=\phi(x_{i})$ (resp., $a=\phi(g)$).
\end{definition}

\begin{lemma}
With the notation as above, let $\tilde{g}$ ($\in A$) be an NB element for each $g\in\Gamma$.

(1) An arbitrary NB element for $g$ is of the form $u\tilde{g}$, where $u\in k^{\times}$.

(2) If $\tilde{x}_{i}$ ($\in A$) is an NB element for $x_{i}$, there exists $c_{g}$ ($\in k$) for each $g$, such that
\begin{equation}
\tilde{g}\tilde{x}_{i}=\chi_{i}(g)\tilde{x}_{i}\tilde{g}+c_{g}\tilde{g}_{i}\tilde{g}.
\end{equation}
\end{lemma}
\begin{proof}
(1) This is easy to see.

(2) This holds since one sees that the element $a:=\tilde{g}\tilde{x}_{i}-\chi_{i}(g)\tilde{x}_{i}\tilde{g}$ is sent to $g_{i}g\otimes a$ by the coaction $A\rightarrow H^{0}\otimes A$.
\end{proof}

Keep $A$, $\tilde{g}$ as above.

\begin{definition}
An NB element $\tilde{x}_{i}$ ($\in A$) for $x_{i}$ is said to be {\em normalized}, if
\begin{equation}
\tilde{g}\tilde{x}_{i}=\chi_{i}(g)\tilde{x}_{i}\tilde{g}\ \ \mbox{for all}\ \ g\in\Gamma.
\end{equation}
By Lemma 3.9.1 this condition is independent of choice of $\tilde{g}$.
\end{definition}

\begin{lemma}
Fix $i\in I$.
Suppose $q_{ii}-1\in k^{\times}$.
Then there exists uniquely a normalized NB element $\tilde{x}_{i}$ ($\in A$) for $x_{i}$.
\end{lemma}
\begin{proof} {\em Existence}.
Choose an NB element $\tilde{x}_{i}$ for $x_{i}$.
By Lemma 3.9.2, we have an element $c$ ($\in k$) such that $\tilde{g}_{i}\tilde{x}_{i}=q_{ii}\tilde{x}_{i}\tilde{g}_{i}+c\tilde{g}_{i}^{2}$.
Set
\[ \hat{x}_{i}=\tilde{x}_{i}+(q_{ii}-1)^{-1}c\tilde{g}_{i}. \]
Then one computes that $\tilde{g}_{i}\hat{x}_{i}=q_{ii}\hat{x}_{i}\tilde{g}_{i}$.
Suppose that a section $\phi$ affords $\tilde{x}_{i}$ so that $\phi(x_{i})=\tilde{x}_{i}$.
We can construct a $k$-linear map $f : S\rightarrow k$ such that
\begin{equation}
f(1)=1,\quad f(x_{j})=\delta_{ij}(q_{ii}-1)^{-1}c\ \ (j\in I).
\end{equation}
Let $\gamma = f\otimes\varepsilon : H^{0}=S\otimes k\Gamma\rightarrow k$, where $\varepsilon$ denotes the counit of $k\Gamma$; this $\gamma$ is convolution-invertible by Lemma 1.14.
The convolution-product $\phi\ast\gamma$ is a section which affords $\hat{x}_{i}$.
We remark that this change of sections does not change the values $\phi(x_{j})$, $\phi(g)$ other than $\phi(x_{i})$.

Fix $g\in\Gamma$.
We have $u\in k^{\times}$ such that $\tilde{g}\tilde{g}_{i}=u\tilde{g}_{i}\tilde{g}$.
Let us write $\tilde{x}_{i}$ for $\hat{x}_{i}$.
Then we have (3.17).
The reduction of $\tilde{g}\tilde{g}_{i}\tilde{x}_{i}$ in two ways shows
\[ q_{ii}c_{g}\tilde{g}_{i}\tilde{g}\tilde{g}_{i}=c_{g}\tilde{g}_{i}\tilde{g}\tilde{g}_{i}, \]
which implies $c_{g}=0$, since $q_{ii}-1\in k^{\times}$.

{\em Uniqueness}. Suppose that $\tilde{x}_{i}$, $\hat{x}_{i}$ are both desired ones.
We see as in the proof of Lemma 3.9.2 that $\tilde{x}_{i}-\hat{x}_{i}=c\tilde{g}_{i}$ for some $c\in k$.
The conjugation by $\tilde{g}_{i}$ gives
\[ (q_{ii}c\tilde{g}_{i}=)\ q_{ii}(\tilde{x}_{i}-\hat{x}_{i})=c\tilde{g}_{i}, \]
which implies $c=0$, since $q_{ii}-1\in k^{\times}$.
\end{proof}

\begin{lemma}
Let $A(\sigma,\mu)$ be the $H^{0}$-cleft object as defined by (3.14).
The images of $x_{i}$ ($i\in I$), $\bar{g}$ ($g\in\Gamma$) under $TV\rtimes_{\sigma}\Gamma\rightarrow A(\sigma,\mu)$ are NB elements in which the images of $x_{i}$ are all normalized.
\end{lemma}
\begin{proof}
This follows by Remark 3.7.
\end{proof}

\begin{definition}
Let
\[ \mathcal{Z}=\mathcal{Z}(\Gamma,\Xi;k) \]
denote the set of all pairs $(\sigma,\mu)$, where $\sigma\in Z^{2}(\Gamma,k^{\times})$, $\mu=(\mu_{ij})\in k^{\Xi(\sigma)}$.
Given two pairs $(\sigma,\mu)$, $(\sigma',\mu')$ in $\mathcal{Z}$, we write $(\sigma,\mu)\sim(\sigma',\mu')$, if there exists a $1$-cochain $\eta : \Gamma\rightarrow k^{\times}$, i.e., a map with $\eta(1)=1$, such that
\begin{eqnarray}
\sigma' & = & \sigma(\partial\eta), \\
\mu_{ij}' & = & \mu_{ij}\eta(g_{i})\eta(g_{j})\quad ((i,j)\in\Xi(\sigma)).
\end{eqnarray}
Recall that (3.20) implies $\Xi(\sigma)=\Xi(\sigma')$.
One sees that $\sim$ defines an equivalence relation on $\mathcal{Z}$.
We denote by
\[ \mathcal{H}(\Gamma,\Xi;k)=\mathcal{Z}/\sim \]
the set of all equivalence classes in $\mathcal{Z}$.
\end{definition}

\begin{proposition}
Suppose that $q_{ii}-1\in k^{\times}$ for all $i\in I$.
Then, $(\sigma,\mu)\mapsto A(\sigma,\mu)$ gives rise to an injection
\[ \mathcal{H}(\Gamma,\Xi;k)\hookrightarrow \Cleft(H^{0}). \]
\end{proposition}
\begin{proof}
Let $(\sigma,\mu)$, $(\sigma',\mu')$ be pairs in $\mathcal{Z}$.
By Lemmas 3.11, 3.12, an isomorphism $A(\sigma',\mu')\xrightarrow{\simeq} A(\sigma,\mu)$ must be induced from an algebra map $f : TV\rtimes_{\sigma'}\Gamma\rightarrow TV\rtimes_{\sigma}\Gamma$ given by
\[ f(x_{i})=x_{i}\ \ (i\in I),\quad f(\bar{g})=\eta(g)\bar{g}\ \ (g\in\Gamma), \]
where $\eta : \Gamma\rightarrow k^{\times}$ is a $1$-cochain.
We see that $f$ indeed induces the desired isomorphism (see Lemma 1.3) if and only if $(\sigma,\mu)\sim(\sigma',\mu')$; this proves the proposition.
In fact, $f$ preserves the relations (3.4) (resp., the relations arising from (3.13)) if and only if (3.20) (resp., (3.21)) holds; it necessarily preserves each element in $\mathcal{R}$.
\end{proof}

\section{The Hopf algebra $H^{\lambda}$ as a cocycle deformation of $H^{0}$}

We fix $\lambda=(\lambda_{ij})\in k^{\Xi}$, throughout this section.

\begin{definition}
Let $\mathcal{R}^{\lambda}$ denote the union of the set $\mathcal{R}$ with the elements
\begin{equation}
(\ad_{c}\,x_{i})(x_{j})-\lambda_{ij}(g_{i}g_{j}-1)\quad ((i,j)\in\Theta)
\end{equation}
in the Hopf algebra $F=TV\dotrtimes\Gamma$ (see (2.2)), where we set $\lambda_{ij}=0$ if $(i,j)\not\in\Xi$.
Define
\[ H^{\lambda}=F/(\mathcal{R}^{\lambda}). \]
\end{definition}

Notice that the element given by (4.1) is a $(1,g_{i}g_{j})$-primitive in $F$. One then sees that $H^{\lambda}$ is a Hopf algebra.
If $\lambda_{ij}=0$ for all $(i,j)\in\Xi$, then $\mathcal{R}^{\lambda}=\mathcal{R}^{0}$, and $H^{\lambda}=H^{0}$; see Definition 2.2.

Recall from Proposition 2.5 that we have the $H^{0}$-cleft object $A(\lambda)$.
Recall from (2.11), (2.12) the definitions of $E$, $e(\lambda)$.

\begin{proposition}
(1) The coproduct $\Delta : F\rightarrow F\otimes F$ induces an algebra map $A(\lambda)\rightarrow A(\lambda)\otimes H^{\lambda}$, by which $A(\lambda)$ is a right $H^{\lambda}$-comodule algebra.

(2) The tensor product of the canonical algebra maps $TV_{r}/(\mathcal{R}_{r})\rightarrow H^{\lambda}$ ($r=1,2,\ldots,N$), composed with the product $(H^{\lambda})^{\otimes N}\rightarrow H^{\lambda}$, induces an isomorphism
\[ e^{\lambda} : E\xrightarrow{\simeq}H^{\lambda} \]
of right $k\Gamma$-module coalgebras.
The composite
\[ \phi^{\lambda}:=e(\lambda)\circ(e^{\lambda})^{-1} : H^{\lambda}\rightarrow A(\lambda) \]
is a right $H^{\lambda}$-colinear section.

(3) $A(\lambda)$ is an $(H^{0},H^{\lambda})$-bicleft object.
\end{proposition}
\begin{proof}
(1) Easy; see the proof of Lemma 2.4.

(2) This together with the fact that $A(\lambda)$ is a right $H^{\lambda}$-cleft object are proved by induction on $L$, just as Proposition 2.5.
The proof in $L=0$ is very similar to the proof of Lemma 2.7.
Suppose that we have proved in $L-1$, and let $\tilde{H}^{\lambda}$, $\tilde{A}(\lambda)$ denote those $H^{\lambda}$, $A(\lambda)$ defined by supposing $\mathcal{R}=\mathcal{R}_{(L-1)}$.
We claim that Lemma 2.8 holds true in $\tilde{H}^{\lambda}$ as well.
In fact we see from the proof of the lemma that the element
\[ z=x_{j}y-q_{ji_{1}}\cdots q_{ji_{l}}yx_{j} \]
is zero in $\tilde{A}(\lambda)$.
But, this element is sent by the right coaction $\tilde{A}(\lambda)\rightarrow\tilde{A}(\lambda)\otimes\tilde{H}^{\lambda}$ to the element
\[ z\otimes 1+gg_{j}\otimes z=gg_{j}\otimes z\ \ \mbox{in}\ \ \tilde{A}(\lambda)\otimes\tilde{H}^{\lambda}, \]
where $g=g_{i_{1}}\cdots g_{i_{l}}$.
This implies that $z=0$ in $\tilde{H}^{\lambda}$, as desired.
Then as in the proof of Proposition 2.5, the induction completes.

(3) It remains to prove that $A(\lambda)$ is an $(H^{0},H^{\lambda})$-bicomodule.
But, this follows since $A(\lambda)$ arises from the $(F,F)$-bicomodule $F$.
\end{proof}

Remark 2.9 can apply to $H^{\lambda}$ as well.

\begin{theorem}
$H^{\lambda}$ and $H^{0}$ are cocycle deformations of each other.
\end{theorem}
\begin{proof}
This follows by Lemma 1.5, since we have the following commutative diagram, in which the unlabeled arrow is the $H^{\lambda}$-coaction on $A(\lambda)$.
\[ \begin{CD}
E @>{\Delta}>> E\otimes E \\
@V{e(\lambda)}VV @VV{e(\lambda)\otimes e^{\lambda}}V \\
A(\lambda) @>>> A(\lambda)\otimes H^{\lambda}
\end{CD} \]
\end{proof}

Recall that $F=TV\dotrtimes\Gamma$ is naturally graded, with $k\Gamma$ its zero component.
Hence it has the induced filtration by $\mathbb{N}=\{0,1,2,\ldots\}$.
The quotient algebra map $F\rightarrow A(\lambda)$ induces a filtration on $A(\lambda)$, with which the map turns into an epimorphism of filtered algebras.
It induces an epimorphism $F\rightarrow \gr A(\lambda)$ of graded algebras.
Similarly, $F\rightarrow H^{\lambda}$ turns into an epimorphism of filtered Hopf algebras, which induces an epimorphism $F\rightarrow \gr H^{\lambda}$ of graded Hopf algebras.
We see that the two graded algebra epimorphisms factor through
\begin{equation}
H^{0}\rightarrow \gr A(\lambda),\quad H^{0}\rightarrow \gr H^{\lambda},
\end{equation}
the latter of which is a graded Hopf algebra map.
Since the $(H^{0},H^{\lambda})$-bicomodule algebra structure maps
\[ H^{0}\otimes A(\lambda)\leftarrow A(\lambda)\rightarrow A(\lambda)\otimes H^{\lambda} \]
on $A(\lambda)$ are filtered, they induce an $(H^{0},\gr H^{\lambda})$-bicomodule algebra structure
\[ H^{0}\otimes \gr A(\lambda) \leftarrow \gr A(\lambda) \rightarrow \gr A(\lambda)\otimes \gr H^{\lambda} \]
on $\gr A(\lambda)$.

\begin{proposition}
(1) $\gr A(\lambda)$ is an $(H^{0},\gr H^{\lambda})$-bicleft object.

(2) The graded algebra map $H^{0}\rightarrow \gr A(\lambda)$ in (4.2) is an isomorphism.

(3) The graded Hopf algebra map $H^{0}\rightarrow \gr H^{\lambda}$ in (4.2) is an isomorphism.
\end{proposition}
\begin{proof}
(1) Notice that $e(\lambda) : E\rightarrow A(\lambda)$ is filtered, whence the section $\phi(\lambda) : H^{0}\rightarrow A(\lambda)$ is a filtered $H^{0}$-colinear map.
It induces an isomorphism $\gr \phi(\lambda) : H^{0}\xrightarrow{\simeq} \gr A(\lambda)$ of graded $H^{0}$-comodules, which is convolution-invertible by Lemma 1.14.
Thus, $\gr A(\lambda)$ is a left $H^{0}$-cleft object with section $\gr \phi(\lambda)$.
Similarly, $\gr A(\lambda)$ is a right $\gr H^{\lambda}$-cleft object with section $\gr \phi^{\lambda}$.

(2) The prescribed map is a map of $H^{0}$-cleft objects, which is necessarily an isomorphism by Lemma 1.3.

(3) It follows by \cite[Theorem 9]{DT} (or see \cite[8.2.4]{Mo}) that $H^{0}$ and $\gr A(\lambda)$, being cleft, are right Galois objects for $H^{0}$, $\gr H^{\lambda}$, respectively.
Combined with the result of Part 2, the (so-called Galois) isomorphisms \cite[8.1.1]{Mo} which ensure the objects to be Galois show that the natural map $H^{0}\otimes H^{0}\rightarrow \gr A(\lambda)\otimes \gr H^{\lambda}$ is an isomorphism.
Since $H^{0}$ ($\simeq \gr A(\lambda)$) includes $k$ as a direct summand, $H^{0}\rightarrow \gr H^{\lambda}$ is an isomorphism.
\end{proof}

Let $\check{\mathcal{R}}(\lambda)$ denote the union of the set $\mathcal{R}$ with the elements
\[ (\ad_{c}\,x_{i})(x_{j})+\lambda_{ij}1\quad ((i,j)\in\Theta) \]
in $F$, where we set $\lambda_{ij}=0$ if $(i,j)\not\in\Xi$.
Define
\[ B(\lambda)=F/(\check{\mathcal{R}}(\lambda)). \]

\begin{proposition}
The obvious $(F,F)$-bicomodule algebra structure on $F$ induces an $(H^{\lambda},H^{0})$-bicomodule algebra structure on $B(\lambda)$.
Moreover, $B(\lambda)$ is an $(H^{\lambda},H^{0})$-bicleft object.
\end{proposition}
\begin{proof}
This can be proved just as Proposition 4.2.
\end{proof}

In the sequel we will write $A$, $B$ for $A(\lambda)$, $B(\lambda)$.

\begin{proposition}
(1) The canonical isomorphism $F\xrightarrow{\simeq} F\Box_{F}F$ induces isomorphisms
\[ \alpha : H^{0}\xrightarrow{\simeq} A\Box_{H^{\lambda}}B,\quad \beta : H^{\lambda}\xrightarrow{\simeq} B\Box_{H^{0}}A \]
of $(H^{0},H^{0})$- and $(H^{\lambda},H^{\lambda})$-bicomodule algebras, respectively.

(2) We have the following commutative diagrams consisting of $(H^{0},H^{\lambda})$- and, resp., $(H^{\lambda},H^{0})$-bicomodule algebras and their isomorphisms, where the unlabeled arrows are canonical isomorphisms.
\begin{center}
\begin{picture}(170,120)
\put(0,105){$H^{0}\Box_{H^{0}}A$}
\put(50,108){\vector(1,0){50}}
\put(65,111){\scriptsize $\alpha\Box\mathrm{id}$}
\put(105,105){$(A\Box_{H^{\lambda}}B)\Box_{H^{0}}A$}
\put(18,52){$A$}
\put(23,67){\vector(0,1){30}}
\put(23,45){\vector(0,-1){30}}
\put(0,0){$A\Box_{H^{\lambda}}H^{\lambda}$}
\put(50,3){\vector(1,0){50}}
\put(65,6){\scriptsize $\mathrm{id}\Box\beta$}
\put(105,0){$A\Box_{H^{\lambda}}(B\Box_{H^{0}}A)$}
\put(130,97){\vector(0,-1){82}}
\end{picture}
\end{center}
\begin{center}
\begin{picture}(170,120)
\put(0,105){$H^{\lambda}\Box_{H^{\lambda}}B$}
\put(50,108){\vector(1,0){50}}
\put(65,111){\scriptsize $\beta\Box\mathrm{id}$}
\put(105,105){$(B\Box_{H^{0}}A)\Box_{H^{\lambda}}B$}
\put(18,52){$B$}
\put(23,67){\vector(0,1){30}}
\put(23,45){\vector(0,-1){30}}
\put(0,0){$B\Box_{H^{0}}H^{0}$}
\put(50,3){\vector(1,0){50}}
\put(65,6){\scriptsize $\mathrm{id}\Box\alpha$}
\put(105,0){$B\Box_{H^{0}}(A\Box_{H^{\lambda}}B)$}
\put(130,97){\vector(0,-1){82}}
\end{picture}
\end{center}

(3) We have the $k$-linear monoidal functors
\begin{eqnarray*}
\Phi_{A}=A\Box_{H^{\lambda}} & : & {}^{H^{\lambda}}\!\!\mathcal{M}\rightarrow{}^{H^{0}}\!\!\mathcal{M}, \\
\Phi_{B}=B\Box_{H^{0}} & : & {}^{H^{0}}\!\!\mathcal{M}\rightarrow{}^{H^{\lambda}}\!\!\mathcal{M},
\end{eqnarray*}
which are mutually quasi-inverses.
Moreover, $\alpha\Box_{H^{0}}$ and $\beta\Box_{H^{\lambda}}$ induce isomorphisms of monoidal functors,
\[ f_{\alpha} : \mathrm{id}\simeq\Phi_{A}\circ\Phi_{B},\quad f_{\beta} : \mathrm{id}\simeq \Phi_{B}\circ\Phi_{A}. \]
$(\Phi_{A},\Phi_{B},f_{\alpha},f_{\beta})$ forms monoidal equivalence data in the sense of Takeuchi \cite{T}; thus, $\Phi_{B}f_{\alpha}=f_{\beta}\Phi_{B}$, $f_{\alpha}\Phi_{A}=\Phi_{A}f_{\beta}$.
\end{proposition}
\begin{proof}
(1) It is easy to see that the bicomodule algebra maps $\alpha$, $\beta$ are induced.
Being maps between cleft objects, they are necessarily isomorphisms, by Lemma 1.3.

(2) This is verified directly.

(3) This follows by Part 2.
\end{proof} 

Let $(\sigma,\mu)\in\mathcal{Z}(\Gamma,\Xi;k)$.
We have the $H^{0}$-cleft object $A(\sigma,\mu)$.
We wish to describe explicitly the corresponding $H^{\lambda}$-cleft object $B\Box_{H^{0}}A(\sigma,\mu)$.
Recall $F=TV\dotrtimes\Gamma$.
Let us write
\[ F(\sigma)=TV\rtimes_{\sigma}\Gamma. \]

\begin{definition}
Let $\mathcal{R}^{\lambda}(\sigma,\mu)$ be the union of the set $\mathcal{R}$ with the elements
\[ (\ad_{c}\,x_{i})(x_{j})+\lambda_{ij}1-\mu_{ij}\bar{g}_{i}\bar{g}_{j}\quad ((i,j)\in\Theta) \]
in $F(\sigma)$, where by convention, $\lambda_{ij}=0$ if $(i,j)\not\in\Xi$, and $\mu_{ij}=0$ if $(i,j)\not\in\Xi(\sigma)$.
Define
\begin{equation}
A^{\lambda}(\sigma,\mu)=F(\sigma)/(\mathcal{R}^{\lambda}(\sigma,\mu)).
\end{equation}
\end{definition}

\begin{proposition}
(1) The natural algebra map $F(\sigma)\rightarrow F\otimes F(\sigma)$ as given in (3.8) induces an algebra map
\[ A^{\lambda}(\sigma,\mu)\rightarrow H^{\lambda}\otimes A^{\lambda}(\sigma,\mu), \]
by which $A^{\lambda}(\sigma,\mu)$ is an $H^{\lambda}$-comodule algebra.
Moreover, $A^{\lambda}(\sigma,\mu)$ is an $H^{\lambda}$-cleft object.

(2) The canonical isomorphism $F(\sigma)\xrightarrow{\simeq}F\Box_{F}F(\sigma)$ induces an isomorphism
\[ A^{\lambda}(\sigma,\mu)\simeq B\Box_{H^{0}}A(\sigma,\mu)=\Phi_{B}(A(\sigma,\mu)) \]
of $H^{\lambda}$-comodule algebras.
\end{proposition}
\begin{proof}
It is easy to prove the first assertion of Part 1.
To prove the remaining, set
\[ C=A(\sigma,\mu),\quad C^{\lambda}=A^{\lambda}(\sigma,\mu). \]
Let $J$ denote the kernel of the composite
\[ F(\sigma)\xrightarrow{\simeq}F\Box_{F}F(\sigma)\rightarrow B\Box_{H^{0}}C \]
of the canonical epimorphisms of $H^{\lambda}$-comodule algebras.
Set $D=F(\sigma)/J$.
We see easily that the quotient map $F(\sigma)\rightarrow D$ factors through an epimorphism
\[ f : C^{\lambda}\rightarrow D \]
of $H^{\lambda}$-comodule algebras.
For Part 2, we should prove that $f$ is an isomorphism.
It suffices to prove that
\[ \Phi_{A}(f) : A\Box_{H^{\lambda}}C^{\lambda}\rightarrow A\Box_{H^{\lambda}}D \]
is an isomorphism.
We see just as above that the isomorphism $F\xrightarrow{\simeq}F\Box_{F}F(\sigma)$ arising from $\mathrm{id} : F\rightarrow F(\sigma)$ (see (3.9)) induces an epimorphism $C\rightarrow A\Box_{H^{\lambda}}C^{\lambda}$.
Its composite with $\Phi_{A}(f)$ coincides with the canonical isomorphism
\[ C\simeq A\Box_{H^{\lambda}}D\ (\simeq (A\Box_{H^{\lambda}}B)\Box_{H^{0}}C), \]
as is easily seen.
It follows that $\Phi_{A}(f)$ is an isomorphism, as desired.
We see from Part 2 and Lemma 1.4 that $A^{\lambda}(\sigma,\mu)$ is an $H^{\lambda}$-cleft object.
\end{proof}

\begin{remark}
We have such a natural section $H^{\lambda}\rightarrow B$ that is defined, just as $\phi^{\lambda}$ in Proposition 4.2.2, as the composite of $(e^{\lambda})^{-1}$ with the natural map $E\rightarrow B$.
Its composite with $B\Box_{H^{0}}\phi(\sigma,\mu)$, where $\phi(\sigma,\mu) : H^{0}\rightarrow A(\sigma,\mu)$ is the section as in (3.15), gives a section,
\begin{equation}
\phi^{\lambda}(\sigma,\mu) : H^{\lambda}\rightarrow A^{\lambda}(\sigma,\mu),
\end{equation}
as is seen from Part 2 above.
Moreover, we see that for each $1\leq r\leq N$, the algebra $E_{r}(\sigma)$ given in (3.16) is naturally embedded into $A^{\lambda}(\sigma,\mu)$, so that $A^{\lambda}(\sigma,\mu)$ together with $\phi^{\lambda}(\sigma,\mu)$ restrict to the $E_{r}$-cleft object $E_{r}(\sigma)$ with the obvious section, in the same sense as in Remark 2.6; see also Remark 3.7.
\end{remark}

\begin{proposition}
Suppose that $q_{ii}-1\in k^{\times}$ for all $i\in I$.
Then, $(\sigma,\mu)\mapsto A^{\lambda}(\sigma,\mu)$ gives rise to an injection
\[ \mathcal{H}(\Gamma,\Xi;k)\hookrightarrow \Cleft(H^{\lambda}). \]
\end{proposition}
\begin{proof}
This follows by Propositions 3.14, 4.8.2.
\end{proof}

\section{Augmented central $H^{\lambda}$-cleft extensions}

We continue to fix $\lambda=(\lambda_{ij})\in k^{\Xi}$, throughout this section.
Define a subset $\Xi_{\lambda}$ of $\Xi$ by
\begin{equation}
\Xi_{\lambda}=\{(i,j)\in\Xi\;|\;\lambda_{ij}\neq 0\}.
\end{equation}
Fix a $k$-module $M$, and regard it as a trivial $\Gamma$-module.
Let $Z^{2}(\Gamma,M)$ denote the additive group of all normalized $2$-cocycles $\Gamma\times\Gamma\rightarrow M$.

\begin{definition}
Let
\[ \mathcal{Z}_{\varepsilon}^{\lambda}=\mathcal{Z}_{\varepsilon}^{\lambda}(\Gamma,\Xi;M) \]
denote the set of all pairs $(s,m)$, where $s\in Z^{2}(\Gamma,M)$, $m=(m_{ij})\in M^{\Xi}$ such that
\begin{equation}
s(g,g_{i})+s(g,g_{j})=s(g_{i},g)+s(g_{j},g)\ \ \mbox{for all}\ \ g\in\Gamma,
\end{equation}
either if $(i,j)\in\Xi_{\lambda}$ or if $m_{ij}\neq 0$; notice that the condition depends on $\lambda$.
Given two pairs $(s,m)$, $(s',m')$ in $\mathcal{Z}_{\varepsilon}^{\lambda}$, we write $(s,m)\sim_{\varepsilon}(s',m')$, if there exists a $1$-cochain $t : \Gamma\rightarrow M$, i.e., a map with $t(1)=0$, such that
\begin{eqnarray*}
s' & = & s+\partial t, \\
m_{ij}' & = & m_{ij}-\lambda_{ij}(t(g_{i})+t(g_{j}))\quad ((i,j)\in\Xi).
\end{eqnarray*}
One sees that $\sim_{\varepsilon}$ defines an equivalence relation on $\mathcal{Z}_{\varepsilon}^{\lambda}$.
We denote by
\[ \mathcal{H}_{\varepsilon}^{\lambda}(\Gamma,\Xi;M)=\mathcal{Z}_{\varepsilon}^{\lambda}/\sim_{\varepsilon} \]
the set of all equivalence classes in $\mathcal{Z}_{\varepsilon}^{\lambda}$.
\end{definition}

Recall from (1.10) the augmented algebra $k_{M}$.
By extending the coefficients from $k$ to $k_{M}$, we define the set $\mathcal{Z}(\Gamma,\Xi;k_{M})$ just as in Definition 3.13.
Let $(\sigma,\mu)\in\mathcal{Z}(\Gamma,\Xi;k_{M})$.
Thus, $\sigma\in Z^{2}(\Gamma,k_{M}^{\times})$, and $\mu=(\mu_{ij})\in (k_{M})^{\Xi(\sigma)}$, where $\Xi(\sigma)$ is defined by (3.12), but $\chi_{i}\chi_{j}=s_{i}^{\sigma}s_{j}^{\sigma}$ should now be regarded as an equation in the group of group homomorphism $\Gamma\rightarrow k_{M}^{\times}$.
We will suppose, as before, $\mu\in (k_{M})^{\Xi}$ by setting $\mu_{ij}=0$ if $(i,j)\not\in\Xi(\sigma)$.
The pair gives rise to the central $H^{\lambda}$-cleft extension $A^{\lambda}(\sigma,\mu)$ over $k_{M}$; see Definition 4.7.
We say that $A^{\lambda}(\sigma,\mu)$ is {\em augmented}, if $x_{i}\mapsto 0$ ($i\in I$), $\bar{g}\mapsto 1$ ($g\in\Gamma$) well defines an algebra map $A^{\lambda}(\sigma,\mu)\rightarrow k$ extending the augmentation $k_{M}\rightarrow k$.

\begin{lemma}
If $(s,m)\in\mathcal{Z}_{\varepsilon}^{\lambda}(\Gamma,\Xi;M)$, the pair $(\sigma,\mu)$ defined by
\[ \sigma=1+s,\quad \mu_{ij}=\lambda_{ij}+m_{ij}\quad ((i,j)\in\Xi) \]
is in $\mathcal{Z}(\Gamma,\Xi;k_{M})$, and $A^{\lambda}(\sigma,\mu)$ is augmented.
Every pair $(\sigma,\mu)$ in $\mathcal{Z}(\Gamma,\Xi;k_{M})$ such that $A^{\lambda}(\sigma,\mu)$ is augmented, thus arises uniquely from a pair in $\mathcal{Z}_{\varepsilon}^{\lambda}(\Gamma,\Xi;M)$.
\end{lemma}
\begin{proof}
One sees easily that such a pair $(\sigma,\mu)$ that satisfies the condition is uniquely of the form as above.
We should prove that the two conditions
\[ \sigma\in Z^{2}(\Gamma,k_{M}^{\times}),\quad \mu_{ij}=0\ \ \mbox{if}\ \ (i,j)\not\in\Xi(\sigma) \]
are equivalent to $(s,m)\in\mathcal{Z}_{\varepsilon}^{\lambda}$.
The first condition is equivalent to $s\in Z^{2}(\Gamma,M)$.
The second is equivalent to that (5.2) holds if $(i,j)\in\Xi_{\lambda}$ or $m_{ij}\neq 0$, since we have
\[ \Xi(\sigma)=\{(i,j)\in\Xi\;|\;\mbox{(5.2) holds}\}. \]
These prove the desired equivalence.
\end{proof}

Given $m=(m_{ij})\in M^{\Xi}$, we define 
\[ \lambda+m=(\lambda_{ij}+m_{ij})\ (\in (k_{M})^{\Xi}). \]

\begin{proposition}
Suppose that $q_{ii}-1\in k^{\times}$ for all $i\in I$.
Then, $(s,m)\mapsto A^{\lambda}(1+s,\lambda+m)$ gives rise to an injection
\[ \mathcal{H}_{\varepsilon}^{\lambda}(\Gamma,\Xi;M)\hookrightarrow \ZCleft_{\varepsilon}(H^{\lambda};k_{M}), \]
where $\ZCleft_{\varepsilon}(H^{\lambda};k_{M})$ is defined by Definition 1.6.
\end{proposition}

For the proof we remark that Definitions 3.8, 3.10 and Lemmas 3.9, 3.11, which are all concerned with cleft extensions for the Hopf algebra $H^{0}$, are generalized in the obvious way to those extensions for $H^{\lambda}$.
In addition, Lemma 3.12 is generalized by the following.

\begin{lemma}
Let $A^{\lambda}(\sigma,\mu)$ be the $H^{\lambda}$-cleft object as defined by (4.3).
The images of $x_{i}$ ($i\in I$), $\bar{g}$ ($g\in\Gamma$) under the quotient map $TV\rtimes_{\sigma}\Gamma\rightarrow A^{\lambda}(\sigma,\mu)$ are NB elements in which the images of $x_{i}$ are all normalized.
\end{lemma}
\begin{proof}
This follows by Remark 4.9.
\end{proof}

\begin{proof}[Proof of Proposition 5.3]
Let $(s,m)$, $(s',m')$ be pairs in $\mathcal{Z}_{\varepsilon}^{\lambda}$.
Set
\[ (\sigma,\mu)=(1+s,\lambda+m),\quad (\sigma',\mu')=(1+s',\lambda+m'). \]
We apply the generalized Lemmas 3.9, 3.11, and Lemma 5.4 above, to $A^{\lambda}(\sigma,\mu)$, $A^{\lambda}(\sigma',\mu')$, regarding them as $H^{\lambda}\otimes k_{M}$-cleft objects.
It then follows that an augmented isomorphism $A^{\lambda}(\sigma',\mu')\xrightarrow{\simeq}A^{\lambda}(\sigma,\mu)$ must be induced from an algebra map $f : TV\rtimes_{\sigma'}\Gamma\rightarrow TV\rtimes_{\sigma}\Gamma$ given by
\[ f(x_{i})=x_{i}\ \ (i\in I),\quad f(\bar{g})=(1+t(g))\bar{g}\ \ (g\in\Gamma), \]
where $t : \Gamma\rightarrow M$ is a $1$-cochain.
We see that $f$ indeed induces the desired isomorphism if and only if $(s,m)\sim_{\varepsilon}(s',m')$; this proves the proposition.
See the proof of Proposition 3.14.
\end{proof}

\section{Classification results in the case of Cartan-type braidings}

We will impose additional assumptions, (C1)--(C4), on $q=(q_{ij})$ and $\mathcal{R}$.
Let $\mathbb{A}=(a_{ij})_{i,j\in I}$ be a GCM (generalized Cartan matrix).
First we assume:
\begin{enumerate}
\renewcommand{\labelenumi}{(C1)}
\item $q$ is {\em of Cartan type} \cite{AS1} with respect to $\mathbb{A}$ in the sense that
\begin{equation}
q_{ij}q_{ji}=q_{ii}^{a_{ij}}\ \ \mbox{for all}\ \ i,j\in I.
\end{equation}
\end{enumerate}

Notice that (6.1) automatically holds if $i=j$, since $a_{ii}=2$.
Since $q_{ij}q_{ji}\neq 1$ implies $a_{ij}\neq 0$, $|i|\neq |j|$ implies that $i$ and $j$ are disconnected in the Dynkin diagram of $\mathbb{A}$.

The following claim follows from the proof of \cite[Lemma A1]{AS1}.

\begin{lemma-definition}
Let $i,j$ be distinct elements in $I$.
If
\begin{equation}
(l)_{q_{ii}}\in k^{\times}\ \ \mbox{for all}\ \ 0<l\leq -a_{ij},
\end{equation}
then the element
\begin{equation}
(\ad_{c}\,x_{i})^{1-a_{ij}}(x_{j})
\end{equation}
in $TV$ is a primitive, where $(l)_{q_{ii}}$ denotes the $q_{ii}$-integer
\[ (l)_{q_{ii}}=1+q_{ii}+\cdots +q_{ii}^{l-1}. \]
{\rm Such a primitive of the form (6.3) is said to be {\em of Serre type}}.
\end{lemma-definition}

Notice that the element (6.3) is of the form (2.5).
In addition to (C1), we assume:
\begin{enumerate}
\renewcommand{\labelenumi}{(C\arabic{enumi})}
\addtocounter{enumi}{1}
\item $\mathcal{R}$ consists only of primitives of Serre type.
More specifically,
\[ \mathcal{R}=\mathcal{R}_{(1)}=\{(\ad_{c}\,x_{i})^{1-a_{ij}}(x_{j})\;|\;(i,j)\in\Pi\}, \]
where $\Pi$ is a subset,
\[ \Pi\subset \{(i,j)\in I\times I\;|\;|i|=|j|,\ i\neq j,\ \mbox{(6.2) holds.}\}, \]
\item $q_{ii}-1\in k^{\times}$ for all $i\in I$,
\item $q_{ii}q_{jj}-q_{ij}q_{ji}\in k^{\times}$ for all $(i,j)\in \Pi$.
\end{enumerate}

Under (C3), the condition (6.2) is equivalent to
\[ q_{ii}^{l}-1\in k^{\times}\ \ \mbox{for all}\ \ 0<l\leq -a_{ij}. \]

\begin{remark}
Let $\sigma\in Z^{2}(\Gamma,k^{\times})$.
We see from (3.3) that
\[ q_{ii}^{\sigma}=q_{ii},\ \ q_{ij}^{\sigma}q_{ji}^{\sigma}=q_{ij}q_{ji}\quad (i,j\in I). \]
Therefore, if $q$, $\mathcal{R}$ satisfy the assumptions (C1)--(C4), then $q^{\sigma}$, ${}_{\sigma}\mathcal{R}$ do.
Cf.\ \cite[Sect.\ 3.2]{AS2}.
\end{remark}

\begin{theorem}
Under the assumptions (C1)--(C4), $(\sigma,\mu)\mapsto A^{\lambda}(\sigma,\mu)$ gives rise to a bijection
\[ \mathcal{H}(\Gamma,\Xi;k)\xrightarrow{\simeq}\Cleft(H^{\lambda}). \]
\end{theorem}
\begin{proof}
By (the proof of) Proposition 4.10, it suffices to prove the surjectivity, supposing $\lambda=0$, or $\lambda_{ij}=0$ for all $(i,j)\in\Xi$.
Let $A$ be an $H^{0}$-cleft object.
We should prove that $A\simeq A(\sigma,\mu)$ for some $(\sigma,\mu)\in\mathcal{Z}$.
Choose a section $\phi : H^{0}\rightarrow A$, and set
\[ \tilde{x}_{i}=\phi(x_{i})\ \ (i\in I),\quad \tilde{g}=\phi(g)\ \ (g\in\Gamma). \]
By Lemma 3.11 and (C3), we may suppose that $\tilde{x}_{i}$ are all normalized.
There uniquely exists $\sigma\in Z^{2}(\Gamma,k^{\times})$ such that
\[ \tilde{g}\tilde{g}'=\sigma(g,g')\widetilde{gg'}\quad (g,g'\in\Gamma). \]
This implies that in ${}_{\check{\sigma}}A$,
\begin{equation}
\tilde{g}\tilde{g}'=\widetilde{gg'}\quad (g,g'\in\Gamma).
\end{equation}
By Proposition 3.4.2 and Remark 6.2, we can replace $V$ with ${}_{\sigma}V$, and may suppose that (6.4) holds in $A$.
It then suffices to prove that $A\simeq A(\mu)$ (see Definition 2.3) for some $\mu\in k^{\Xi}$.
We see that $x_{i}\mapsto \tilde{x}_{i}$ ($i\in I$), $g\mapsto\tilde{g}$ ($g\in\Gamma$) define an $H^{0}$-comodule algebra map $f : F=TV\dotrtimes\Gamma\rightarrow A$.
It suffices to prove that $f$ factors through $A(\mu)\rightarrow A$, which will be necessarily an isomorphism by Lemma 1.3.
Since each $(\ad_{c}\,x_{i})(x_{j})$, where $(i,j)\in\Theta$, is a $(1,g_{i}g_{j})$-primitive in $F$, we have $\mu=(\mu_{ij})\in k^{\Theta}$ such that
\begin{equation}
\tilde{x}_{i}\tilde{x}_{j}=q_{ij}\tilde{x}_{j}\tilde{x}_{i}+\mu_{ij}\tilde{g}_{i}\tilde{g}_{j}\quad ((i,j)\in\Theta).
\end{equation}
Here, $\mu_{ij}=0$ if $(i,j)\not\in\Xi$; otherwise, $\mu_{ij}\tilde{g}_{i}\tilde{g}_{j}$ ($=\tilde{x}_{i}\tilde{x}_{j}-q_{ij}\tilde{x}_{j}\tilde{x}_{i}$) would not then commute with all $\tilde{g}$, where $g\in\Gamma$.
Choose an arbitrary element (6.3) from $\mathcal{R}$, and let $y$ denote its image in $A$.
Then we have as above, $c\in k$ such that
\[ y=c\tilde{g}_{i}^{1-a_{ij}}\tilde{g}_{j}. \]
The conjugations by $\tilde{g}_{i}$, and by $\tilde{g}_{j}$ give
\begin{eqnarray}
q_{ii}^{1-a_{ij}}q_{ij}y & = & c\tilde{g}_{i}^{1-a_{ij}}\tilde{g}_{j}\ (=y), \\
q_{ji}^{1-a_{ij}}q_{jj}y & = & c\tilde{g}_{i}^{1-a_{ij}}\tilde{g}_{j}\ (=y).
\end{eqnarray}
By (6.1), (6.6), one sees that
\begin{equation}
\frac{q_{ii}}{q_{ji}}y=y.
\end{equation}
This together with (6.1), (6.7) imply that
\begin{equation}
\frac{q_{jj}}{q_{ij}}y=\frac{q_{jj}}{q_{ij}}\left(\frac{q_{ii}}{q_{ji}}\right)^{a_{ij}}y=q_{ji}^{1-a_{ij}}q_{jj}y=y.
\end{equation}
Since one sees by (6.8), (6.9) that $(q_{ii}q_{jj}-q_{ij}q_{ji})y=0$, it follows by (C4) that $y=0$.
We have shown that each element in $\mathcal{R}$ is zero in $A$.
This together with (6.5) prove the desired factorization.
\end{proof}

\begin{theorem}
Under the assumptions (C1)--(C4), $(s,m)\mapsto A^{\lambda}(1+s,\lambda+m)$ gives rise to a bijection
\[ \mathcal{H}_{\varepsilon}^{\lambda}(\Gamma,\Xi;M)\xrightarrow{\simeq}\ZCleft_{\varepsilon}(H^{\lambda};k_{M})\ (\simeq H^{2}(H^{\lambda},M)). \]
\end{theorem}
\begin{proof}
By Proposition 5.3, it suffices to prove the surjectivity.
Let $A=(A,\varepsilon_{A})$ be an augmented central $H^{\lambda}$-cleft extension over $k_{M}$.
By Lemma 1.7, it has an augmented section $\phi : H^{\lambda}\rightarrow A$.
Set $\tilde{x}_{i}=\phi(x_{i})$ ($i\in I$), $\tilde{g}=\phi(g)$ ($g\in\Gamma$).
Then,
\begin{equation}
\varepsilon_{A}(\tilde{x}_{i})=0,\ \ \varepsilon_{A}(\tilde{g})=1\quad (i\in I,\ g\in\Gamma).
\end{equation}
Regard $A$ as an $H^{\lambda}\otimes k_{M}$-cleft object.
Then, $\tilde{x}_{i}$, $\tilde{g}$ form NB elements.
By (C3), the $\tilde{x}_{i}$ can be so chosen as to be normalized.
This is possible with the property (6.10) preserved, as is seen from the proof of Lemma 3.11.
Notice that Theorem 6.3 still holds true over $k_{M}$, under the same assumptions.
Combining the result with Lemmas 3.9.1, 3.11 (generalized to $H^{\lambda}$), we have $(\sigma,\mu)\in\mathcal{Z}(\Gamma,\Xi;k_{M})$ and $1$-cochain $\eta : \Gamma\rightarrow k_{M}^{\times}$, such that $x_{i}\mapsto\tilde{x}_{i}$ ($i\in I$), $g\mapsto\eta(g)^{-1}\tilde{g}$ ($g\in\Gamma$) induce an isomorphism $A^{\lambda}(\sigma,\mu)\simeq A$ of central $H^{\lambda}$-cleft extensions over $k_{M}$.
By replacing $(\sigma,\mu)$ with $(\sigma',\mu')$ defined by (3.20), (3.21), we may suppose that $\eta(g)=1$ for all $g\in\Gamma$; see the proof of Proposition 3.14.
Then, $A^{\lambda}(\sigma,\mu)$, being augmented, is necessarily of the form $A^{\lambda}(1+s,\lambda+m)$, by Lemma 5.2.
This proves the desired surjectivity.
\end{proof}

\section{Results for the quantized enveloping algebra $U_{q}$ and the Borel subalgebra $B_{q}$}

In this section we suppose that $k$ is a field. Let $\mathbb{A}=(a_{ij})_{1\leq i,j\leq n}$ be an $n\times n$ GCM which is symmetrized by a diagonal matrix $\mathbb{D}=\diag(\cdots d_{i}\cdots)$ with $0<d_{i}\in\mathbb{Z}$, so that $d_{i}a_{ij}=d_{j}a_{ji}$ for all $1\leq i,j\leq n$.
Let $(\mathfrak{h},\Pi,\Pi^{\vee})$ be a realization of $\mathbb{A}$ (defined over $\mathbb{C}$), with
\[ \Pi=\{\alpha_{1},\dotsc,\alpha_{n}\}\ (\subset\mathfrak{h}^{\ast}),\quad \Pi^{\vee}=\{\alpha_{1}^{\vee},\dotsc,\alpha_{n}^{\vee}\}\ (\subset\mathfrak{h}). \]
Let $\mathfrak{h}_{\mathbb{Z}}\subset\mathfrak{h}$ be a lattice such that
\[ d_{i}\alpha_{i}^{\vee}\in\mathfrak{h}_{\mathbb{Z}},\quad \alpha_{i}(\mathfrak{h}_{\mathbb{Z}})\subset\mathbb{Z}\quad (1\leq i\leq n). \]
Fix $0\neq q\in k$.
Set
\[ q_{i}=q^{d_{i}}\quad (1\leq i\leq n), \]
and suppose that for each $1\leq i\leq n$,
\begin{equation}
\ord(q_{i}^{2})>\max\{1,-a_{ij}\;|\;j\neq i\}.
\end{equation}
Here, $\ord(q_{i}^{2})$ denotes the order of $q_{i}^{2}$ in the group $k^{\times}$;
it is infinite (and then satisfies (7.1)), if $q_{i}$ is not a root of $1$.

\begin{definition}{\it \cite{Ji}}
The {\em quantized enveloping algebra} $U_{q}$ is the algebra (over $k$) which is generated by the elements
\[ K^{h}\ (h\in\mathfrak{h}_{\mathbb{Z}}),\ \ X_{i}^{+},X_{i}^{-}\ (1\leq i\leq n), \]
and is defined by the relations
\begin{eqnarray}
K^{h}K^{h'} & = & K^{h+h'}\quad (h,h'\in\mathfrak{h}_{\mathbb{Z}}), \nonumber \\
K^{0} & = & 1, \nonumber \\
K^{h}X_{i}^{\pm}K^{-h} & = & q^{\pm \alpha_{i}(h)}X_{i}^{\pm}, \\
X_{i}^{+}X_{j}^{-}-X_{j}^{-}X_{i}^{+} & = & \delta_{ij}\frac{g_{i}-g_{i}^{-1}}{q_{i}-q_{i}^{-1}},
\end{eqnarray}
\begin{equation}
\sum_{r=0}^{1-a_{ij}}(-1)^{r}\left[\begin{array}{c} 1-a_{ij} \\ r \end{array}\right]_{i}(X_{i}^{\pm})^{1-a_{ij}-r}X_{j}^{\pm}(X_{i}^{\pm})^{r}=0\quad (i\neq j).
\end{equation}
Here, the $g_{i}^{\pm 1}$ in (7.3) are defined by
\[ g_{i}^{\pm 1}=K^{\pm d_{i}\alpha_{i}^{\vee}}\quad (1\leq i\leq n), \]
and the symbol $\left[\begin{array}{c} \ \ \\ \ \ \end{array}\right]_{i}$ denotes the Gauss binomial coefficients defined by
\begin{eqnarray*}
\left[\begin{array}{c} m \\ r \end{array}\right]_{i} & = & \frac{[m]_{i}[m-1]_{i}\cdots[m-r+1]_{i}}{[r]_{i}[r-1]_{i}\cdots [1]_{i}}\quad (m\geq r > 0), \\
\left[\begin{array}{c} m \\ 0 \end{array}\right]_{i} & = & 1,\ \ \mbox{where}\ \ [m]_{i}=\frac{q_{i}^{m}-q_{i}^{-m}}{q_{i}-q_{i}^{-1}}\quad (m>0).
\end{eqnarray*}
The {\em Borel subalgebra} $B_{q}$ of $U_{q}$ is the subalgebra generated by $K^{h}$ ($h\in\mathfrak{h}_{\mathbb{Z}}$), $X_{i}^{+}$ ($1\leq i\leq n$).
\end{definition}

The assumption (7.1) ensures that the denominators of the right-hand side of (7.3), and of the Gauss binomial coefficients in (7.4) are non-zero.
It is known that $U_{q}$ forms a pointed Hopf algebra including $B_{q}$ as a Hopf subalgebra, in which $K^{h}$ are grouplikes, $X_{i}^{+}$ is a $(1,g_{i})$-primitive, and $X_{i}^{-}$ is a $(g_{i}^{-1},1)$-primitive.
We will see below that $B_{q}$ (resp.\, $U_{q}$) is of the form $H^{0}$ (resp.\, $H^{\lambda}$).

Let $m=\dim_{\mathbb{C}}\mathfrak{h}$.
Then, $m=2n-\rank\mathbb{A}$, and $\mathfrak{h}_{\mathbb{Z}}\simeq\mathbb{Z}^{m}$.
Define
\begin{equation}
\Gamma := \{K^{h}\in U_{q}\;|\;h\in\mathfrak{h}_{\mathbb{Z}}\}.
\end{equation}
Then, $\Gamma$ is a free abelian group of rank $m$, such that $\mathfrak{h}_{\mathbb{Z}}\simeq\Gamma$ via $h\mapsto K^{h}$.
Moreover, it coincides with the groups $G(U_{q}),G(B_{q})$ of the grouplikes in the Hopf algebras $U_{q},B_{q}$.
Since $\alpha_{j}(\alpha_{i}^{\vee})=a_{ij}$ by definition, we have by (7.2),
\begin{equation}
g_{i}X_{j}^{\pm}=q_{i}^{a_{ij}}X_{j}^{\pm}g_{i}.
\end{equation}

\begin{notation}
Let $I_{+}=\{1,2,\ldots,n\}$.
Define $\chi_{i}\in\hat{\Gamma}$ ($i\in I_{+}$) by
\[ \chi_{i}(K^{h})=q^{\alpha_{i}(h)}\quad (h\in\mathfrak{h}_{\mathbb{Z}}). \]
Define an object $V_{+}$ in ${}^{\Gamma}_{\Gamma}\YD$ by
\begin{equation}
V_{+}=\{x_{i},g_{i},\chi_{i}\}_{i\in I_{+}}.
\end{equation}
Since $\chi_{j}(g_{i})=q_{i}^{a_{ij}}$, the associated braiding is given by
\begin{equation}
c : V_{+}\otimes V_{+} \xrightarrow{\simeq} V_{+}\otimes V_{+},\quad c(x_{i}\otimes x_{j})=q_{i}^{a_{ij}}x_{j}\otimes x_{i},
\end{equation}
which is afforded by the matrix
\begin{equation}
\mathbf{q}_{+}=(q_{ij})_{i,j\in I_{+}},\ \ \mbox{where}\ \ q_{ij}=q_{i}^{a_{ij}}=q^{d_{i}a_{ij}}.
\end{equation}
\end{notation}

We see that $\mathbf{q}_{+}$ is of Cartan type with respect to the GCM $\mathbb{A}$; see (C1) in Section 6.
To construct, as in Section 2, $H^{0}$ from the $V_{+}$ above, we do not decompose $I_{+}$ into disconnected subsets, so that we set $\mathcal{R}^{0}=\mathcal{R}$.
We choose as $\mathcal{R}$ ($=\mathcal{R}_{(1)}$) the set
\begin{equation}
\mathcal{R}_{+}=\{(\ad_{c}\,x_{i})^{1-a_{ij}}(x_{j})\;|\;i\neq j\},
\end{equation}
which consists of primitives in $TV_{+}$ of Serre type, by the assumption (7.1); see Lemma 6.1.
We thus construct
\[ H^{0}=TV_{+}\dotrtimes\Gamma/(\mathcal{R}_{+}). \]

\begin{lemma}
$x_{i}\mapsto X_{i}^{+}$ ($i\in I_{+}$) gives such an isomorphism $H^{0}\xrightarrow{\simeq} B_{q}$ of Hopf algebras that is identical on the group algebra $k\Gamma$.
\end{lemma}
\begin{proof}
This is easy to see, and is well known.
Notice that the image of $(\ad_{c}\,x_{i})^{1-a_{ij}}(x_{j})$ in $B_{q}$ coincides with the left-hand side of (7.6) with plus-sign.
\end{proof}

In some places we will impose the following assumption.
\begin{enumerate}
\renewcommand{\labelenumi}{($\mathrm{C}_{q}$)}
\item $q^{2(d_{i}+d_{j}-d_{i}a_{ij})}\neq 1$ for all $i\neq j$ in $I_{+}$
\end{enumerate}
\begin{theorem}
Under the assumption ($\mathrm{C}_{q}$), the restriction maps give a bijection
\[ \ZCleft(B_{q};Z) \xrightarrow{\simeq} \Cleft(Z\Gamma)\ (\simeq (Z^{\times})^{\mbox{\rm\normalsize (\raisebox{-0.5ex}{$\stackrel{\scriptstyle m}{\scriptstyle 2}$})}}), \]
and a $k$-linear isomorphism
\[ H^{2}(B_{q},M)\xrightarrow{\simeq}H^{2}(\Gamma,M)\ (\simeq M^{\mbox{\rm\normalsize (\raisebox{-0.5ex}{$\stackrel{\scriptstyle m}{\scriptstyle 2}$})}}), \]
where $Z\neq 0$ is a commutative algebra, and $M$ is an arbitrary $k$-vector space, regarded as the trivial module over $B_{q}$, and over $\Gamma$.
\end{theorem}
\begin{proof}
In the present situation,
\[ \mathcal{H}(\Gamma,\Xi;Z)=H^{2}(\Gamma,Z^{\times}),\quad \mathcal{H}_{\varepsilon}^{\lambda}(\Gamma,\Xi;M)=H^{2}(\Gamma,M), \]
since $\Theta=\Xi=\emptyset$.
The theorem will follow from Theorems 6.3, 6.4, if we see that the assumptions (C1)--(C4) are satisfied.
Obviously, (C1), (C2) are satisfied.
(C3) is ensured by (7.1).
Since we see from (7.9) that
\[ \frac{q_{ii}q_{jj}}{q_{ij}q_{ji}}=q^{2(d_{i}+d_{j}-d_{i}a_{ij})}, \]
(C4) is precisely ($\mathrm{C}_{q}$).
\end{proof}

\begin{remark}
Let $M$ be a trivial module over $\Gamma$ ($\simeq \mathbb{Z}^{m}$) as above.
Choose an ordered $\mathbb{Z}$-free basis of $\Gamma$, say $e_{1},\dotsc,e_{m}$.
As is well known,
\[ \sigma\mapsto (\sigma(e_{j},e_{i})-\sigma(e_{i},e_{j}))_{i<j} \]
gives an isomorphism
\begin{equation}
H^{2}(\Gamma,M)\simeq M^{\mbox{\normalsize (\raisebox{-0.5ex}{$\stackrel{\scriptstyle m}{\scriptstyle 2}$})}}.
\end{equation}
In particular we have
\[ \Cleft(Z\Gamma)\simeq H^{2}(\Gamma,Z^{\times})\simeq (Z^{\times})^{\mbox{\normalsize (\raisebox{-0.5ex}{$\stackrel{\scriptstyle m}{\scriptstyle 2}$})}}, \]
where $Z\neq 0$ is a commutative algebra.
\end{remark}

\begin{notation}
Let $I_{-}=-I_{+}=\{-1,-2,\ldots,-n\}$.
Define $V_{-}\in {}^{\Gamma}_{\Gamma}\YD$ by
\[ V_{-}=\{x_{-i},g_{i},\chi_{i}^{-1}\}_{i\in I_{+}}. \]
Set
\[ I=I_{+}\sqcup I_{-},\quad V=V_{+}\oplus V_{-}\ (\in{}^{\Gamma}_{\Gamma}\YD). \]
We denote the associated braiding on $V$, which extends the $c$ on $V_{+}$ given by (7.8), by the same symbol $c : V\otimes V\xrightarrow{\simeq} V\otimes V$.
\end{notation}

Notice that $(V,c)$ is of diagonal type, with $c$ afforded by the matrix
\[ \mathbf{q}=\left(\begin{array}{cc} \mathbf{q}_{+} & \check{\mathbf{q}}_{+} \\ \mathbf{q}_{+} & \check{\mathbf{q}}_{+} \end{array}\right),\ \ \mbox{where}\ \ \check{\mathbf{q}}_{+}=(q_{ij}^{-1})_{i,j\in I_{+}}. \]
Here we have supposed that the indices in $I$ are in the order $1,2,\ldots,n,-1,-2,\ldots,-n$.
Moreover, $\mathbf{q}$ is of Cartan type with respect to the GCM $\mathbb{A}\oplus\mathbb{A}$.

To construct, as in Section 4, $H^{\lambda}$ from $V$, we decompose $I$ into $I_{1}=I_{-}$ and $I_{2}=I_{+}$, which are mutually disconnected.
Then,
\[ \Theta = \{(i,-j)\;|\;i,j\in I_{+}\},\quad \Xi\supset \{(i,-i)\;|\; i\in I_{+}\}. \]
We can choose $\lambda=(\lambda_{i,-j})\in k^{\Xi}$, so as 
\begin{equation}
\lambda_{i,-j}=\delta_{ij}\frac{1}{q_{i}-q_{i}^{-1}}\quad (i,j\in I_{+}).
\end{equation}
Define 
\[ \mathcal{R}=\mathcal{R}_{(1)}=\{(\ad_{c}\,x_{\pm i})^{1-a_{ij}}(x_{\pm j})\;|\;i,j\in I_{+},i\neq j\}, \]
which consists of primitives in $TV$ of Serre type, by the assumption (7.1).
Therefore, $\mathcal{R}^{\lambda}$ denotes the union of the set $\mathcal{R}$ with the elements
\[ (\ad_{c}\,x_{i})(x_{-j})-\delta_{ij}\frac{g_{i}g_{j}-1}{q_{i}-q_{i}^{-1}}\quad (i,j\in I_{+}). \]
We thus construct
\[ H^{\lambda} = TV\dotrtimes\Gamma/(\mathcal{R}^{\lambda}). \]

\begin{lemma} $x_{i}\mapsto X_{i}^{+}$, $x_{-i}\mapsto X_{i}^{-}K_{i}$ ($i\in I_{+}$) give such an isomorphism $H^{\lambda}\stackrel{\simeq}{\rightarrow}U_{q}$ of Hopf algebras that is identical on $k\Gamma$.
\end{lemma}
\begin{proof}
This is also well known; a detailed proof can be found in \cite[Bemerkung 4.2]{G1}.
We only remark that the image of $(\ad_{c}\,x_{-i})^{1-a_{ij}}(x_{-j})$ coincides with the left-hand side of (7.4) with minus-sign, modulo multiplication by a unit.
\end{proof}

The isomorphism above transfers the natural filtration on $H^{\lambda}$ onto $U_{q}$, so that $U_{q}$ turns into a filtered Hopf algebra with respect to the filtration $(U_{q})_{0}\subset (U_{q})_{1}\subset\cdots$ in which $(U_{q})_{r}$ is spanned by those monomials in $K^{h}$ ($h\in\mathfrak{h}_{\mathbb{Z}}$), $X_{i}^{\pm}$ ($i\in I_{+}$) containing $X_{i}^{\pm}$ at most $r$.
Let $\gr U_{q}$ denote the associated graded Hopf algebra.
On the other hand, let $(U_{q})^{0}$ denote the (graded) Hopf algebra defined just as $U_{q}$ except that the relation (7.3) is replaced by
\[ X_{i}^{+}X_{j}^{-}-X_{j}^{-}X_{i}^{+}=0. \]

\begin{theorem}
(1) We have a canonical isomorphism $(U_{q})^{0}\simeq\gr U_{q}$ of graded Hopf algebras.

(2) $U_{q}$ and $(U_{q})^{0}$ are cocycle deformations of each other.
\end{theorem}
\begin{proof}
As in Section 2, we define $H^{0}=TV\dotrtimes\Gamma/(\mathcal{R}^{0})$, where $\mathcal{R}^{0}$ is the union of the set $\mathcal{R}$ with the elements $(\ad_{c}\,x_{i})(x_{-j})$ ($i,j\in I_{+}$).
Just as Lemma 7.7, we have an isomorphism $H^{0}\simeq (U_{q})^{0}$ of graded Hopf algebras.
The theorem now follows by Theorem 4.3 and Proposition 4.4.3.
\end{proof}

\begin{remark}
Suppose that the GCM $\mathbb{A}$ is of finite type.
In this case the result of Part 2 above was proved by Kassel and Schneider \cite{KS} in a different method.
They also proved the same result for the Frobenius-Lusztig kernels $u_{q}$, which is a finite-dimensional quotient Hopf algebra of $U_{q}$ defined when $q$ is a root of $1$;
this result also follows from our Theorem 4.3, in the same way as above.
Didt \cite{Di} obtained results on cocycle deformations of finite-dimensional Nichols algebras and their liftings.
Among others, he obtained the result \cite[Theorem 1]{Di} for those algebras defined with zero {\em root vector parameters} \cite{AS2,AS3}, which generalizes the above-cited result for $u_{q}$.
Even this result by Didt follows from our Theorem 4.3; see Remark A.2 in Appendix.
\end{remark}

\begin{definition}
Let $U_{q}'$ denote the subalgebra of $U_{q}$ generated by $g_{i}^{\pm 1}$, $X_{i}^{\pm}$ ($i\in I_{+}$);
this is indeed a Hopf subalgebra.
\end{definition}

If $\det\mathbb{A}\neq 0$ (e.g., $\mathbb{A}$ is of finite type), then $\mathfrak{h}_{\mathbb{Z}}$ can be chosen so as $\mathfrak{h}_{\mathbb{Z}}=\bigoplus_{i=1}^{n}\mathbb{Z}d_{i}\alpha_{i}^{\vee}$, in which case $U_{q}'=U_{q}$.

Replace the definition (7.5) of $\Gamma$ by
\[ \Gamma := \{g_{i}^{l}\in U_{q}'\;|\;l\in\mathbb{Z},\ i\in I_{+}\}. \]
This is the free abelian group on the set $\{g_{i}\}_{i\in I_{+}}$, whence $\Gamma\simeq\mathbb{Z}^{n}$, and it coincides with the group $G(U_{q}')$ of the grouplikes in $U_{q}'$.
Accordingly, replacing the definition of $H^{\lambda}$, let $H^{\lambda}$ denote the Hopf subalgebra of the original one which is generated by $g_{i}^{\pm 1}$, $x_{\pm i}$ ($i\in I_{+}$).
We then still have just as in Lemma 7.7, such a Hopf algebra isomorphism $H^{\lambda}\simeq U_{q}'$ that is identical on $k\Gamma$.

Let $Z\neq 0$ be a commutative algebra.
Define
\[ \nabla = \{(i,j)\;|\;i,j\in I_{+},\ i<j\}. \]
Then, $\#\nabla=\left(\begin{array}{c} n \\ 2 \end{array}\right)$.
Let $u=(u_{ij})_{i<j}\in (Z^{\times})^{\nabla}$.
We set
\[ u_{ii}=1,\ \ u_{ji}=u_{ij}^{-1}\ \ (i<j). \]
Define
\begin{equation}
\Xi(u)=\{(i,j)\in I_{+}^{2}\;|\; q^{d_{r}(a_{ri}-a_{rj})}=u_{ir}u_{jr}\ \ \mbox{for all}\ \ r\in I_{+}\}.
\end{equation}
The group $(Z^{\times})^{I_{+}}$ acts on $Z^{\Xi(u)}$ from the right by
\[ (\mu^{\omega})_{ij}=\mu_{ij}\omega_{i}\omega_{j}\quad ((i,j)\in\Xi(u)), \]
where $\omega = (\omega_{i})\in (Z^{\times})^{I_{+}}$, $\mu=(\mu_{ij})\in Z^{\Xi(u)}$.
Let $Z^{\Xi(u)}/(Z^{\times})^{I_{+}}$ denote the set of all orbits.

\begin{example}
Assume ($\mathrm{C}_{q}$).
If $u=1$, or namely if $u_{ij}=1$ for all $i<j$, then
\[ \Xi(1)=\{(i,i)\;|\;i\in I_{+}\}, \]
since the condition given in (7.13) for $r=i,j$ implies $q^{2(d_{i}+d_{j}-d_{i}a_{ij})}=1$.
It follows that
\[ Z^{\Xi(1)}/(Z^{\times})^{I_{+}}\simeq (Z/(Z^{\times})^{2})^{n}. \]
\end{example}

\begin{definition}
Given $\mu=(\mu_{ij})\in Z^{\Xi(u)}$, let $A_{q}(u,\mu)$ denote the $Z$-algebra which is generated by the elements $\tilde{g}_{i}^{\pm 1},\tilde{X}_{i}^{\pm}$ ($i\in I_{+}$), and is defined by the same relation as (7.6) (with $g_{i}^{\pm 1},X_{i}^{\pm}$ replaced by $\tilde{g}_{i}^{\pm 1}, \tilde{X}_{i}^{\pm}$), and by the additional relations
\[ \tilde{g}_{i}\tilde{g}_{i}^{-1}=1=\tilde{g}_{i}^{-1}\tilde{g}_{i}, \]
\begin{eqnarray*}
\tilde{g}_{j}\tilde{g}_{i} & = & u_{ij}\tilde{g}_{i}\tilde{g}_{j}\quad (i<j), \\
\tilde{X}_{i}^{+}\tilde{X}_{j}^{-}-\tilde{X}_{j}^{-}\tilde{X}_{i}^{+} & = & \mu_{ij}\tilde{g}_{i}-\delta_{ij}\frac{\tilde{g}_{i}^{-1}}{q_{i}-q_{i}^{-1}},
\end{eqnarray*}
\begin{eqnarray}
\sum_{r=0}^{1-a_{ij}}(-1)^{r}\left[\begin{array}{c} 1-a_{ij} \\ r \end{array}\right]_{i}(\tilde{X}_{i}^{+})^{1-a_{ij}-r}\tilde{X}_{j}^{+}(\tilde{X}_{i}^{+})^{r} & = & 0\ (i\neq j), \nonumber \\
\sum_{r=0}^{1-a_{ij}}(-u_{ij})^{r}\left[\begin{array}{c} 1-a_{ij} \\ r \end{array}\right]_{i}(\tilde{X}_{i}^{-})^{1-a_{ij}-r}\tilde{X}_{j}^{-}(\tilde{X}_{i}^{-})^{r} & = & 0\ (i\neq j).
\end{eqnarray}
\end{definition}

\begin{proposition}
$\rho(\tilde{g}_{i}^{\pm 1})=g_{i}^{\pm 1}\otimes \tilde{g}_{i}^{\pm 1}$, $\rho(\tilde{X}_{i}^{+})=X_{i}^{+}\otimes 1+g_{i}\otimes\tilde{X}_{i}^{+}$, $\rho(\tilde{X}_{i}^{-})=X_{i}^{-}\otimes\tilde{g}_{i}^{-1}+1\otimes\tilde{X}_{i}^{-}$ define a $Z$-algebra map
\[ \rho : A_{q}(u,\mu)\rightarrow U_{q}'\otimes A_{q}(u,\mu), \]
by which $A_{q}(u,\mu)$ is a central $U_{q}$-cleft extensions over $Z$.
\end{proposition}
\begin{proof}
Choose $\sigma\in Z^{2}(\Gamma,Z^{\times})$ whose cohomology class corresponds to $u$ via a bijection as in (7.11);
thus, $\sigma(g_{j},g_{i})/\sigma(g_{i},g_{j})=u_{ij}$.
Then, $\Xi(u)\simeq\Xi(\sigma)$ via $(i,j)\mapsto (i,-j)$, by which we can regard $\mu\in \Xi(\sigma)$.
Then one can construct the $H^{\lambda}\otimes Z$-cleft object $A^{\lambda}(\sigma,\mu)$ as in Definition 4.7, where $\lambda$ is kept as in (7.12).
Notice that (7.14) is equivalent to
\[ \sum_{r=0}^{1-a_{ij}}(-1)^{r}\left[\begin{array}{c} 1-a_{ij} \\ r \end{array}\right]_{i}(\tilde{X}_{i}^{-}\tilde{K}_{i})^{1-a_{ij}-r}(\tilde{X}_{j}^{-}\tilde{K}_{j})(\tilde{X}_{i}^{-}\tilde{K}_{i})^{r}=0. \]
Then one can see that the $Z$-algebra map $T(V\otimes Z)\rtimes_{\sigma}\Gamma\rightarrow A_{q}(u,\mu)$ given by $\bar{g}_{i}^{\pm 1}\mapsto \tilde{g}_{i}^{\pm 1}$, $x_{i}\mapsto \tilde{X}_{i}^{+}$, $x_{-i}\mapsto \tilde{X}_{i}^{-}\tilde{g}_{i}$ ($i\in I_{+}$) induces a $Z$-algebra isomorphism $A^{\lambda}(\sigma,\mu)\simeq A_{q}(u,\mu)$, by which and by $H^{\lambda}\simeq U_{q}'$, the $H^{\lambda}$-comodule algebra structure on $A^{\lambda}(\sigma,\mu)$ is transferred to the $U_{q}'$-comodule algebra structure on $A_{q}(u,\mu)$ prescribed above.
This proves the proposition.
\end{proof}

\begin{theorem}
Assume ($\mathrm{C}_{q}$).
Then the maps $Z^{\Xi(u)}\rightarrow\ZCleft(U_{q}';Z)$ ($u\in (Z^{\times})^{\nabla}$) defined by $\mu\mapsto A_{q}(u,\mu)$ give rise to a bijection
\[ \bigsqcup_{u\in (Z^{\times})^{\nabla}}Z^{\Xi(u)}/(Z^{\times})^{I_{+}}\xrightarrow{\simeq}\ZCleft(U_{q}';Z). \]
\end{theorem}
\begin{proof}
The assumptions (C1)--(C4) are satisfied, as is seen in a similar way of the proof of Theorem 7.4.
We remark in particular that (C4) is still equivalent to ($\mathrm{C}_{q}$) in the present situation.

Given $(\sigma,\mu)\in\mathcal{Z}(\Gamma,\Xi;Z)$, suppose that the cohomology class of $\sigma$ corresponds to $u\in (Z^{\times})^{\nabla}$. Assign then $\mu\in Z^{\Xi(u)}$ to $(\sigma,\mu)$, identifying $\Xi(u)=\Xi(\sigma)$.
One sees that this assignment gives rise to a map
\[ \mathcal{H}(\Gamma,\Xi;Z)\rightarrow \bigsqcup_{u}Z^{\Xi(u)}/(Z^{\times})^{\nabla}. \]
By Theorem 6.3 and the proof of Proposition 7.13, it suffices to prove that this map is bijective.
The surjectivity is obvious.
For the injectivity, suppose that $(\sigma,\mu)$, $(\sigma',\mu')$ have the same image.
Then, $\sigma$ and $\sigma'$ are cohomologous, and we have $\eta(g_{i})\in Z^{\times}$ ($i\in I_{+}$) such that
\[ \mu_{ij}'=\mu_{ij}\eta(g_{i})\eta(g_{j})\quad ((i,j)\in\Xi(u)). \]
We should prove that $g_{i}\mapsto \eta(g_{i})$ can extend to a $1$-cochain $\eta : \Gamma\rightarrow Z^{\times}$ such that $\sigma'=\sigma(\partial\eta)$.
This is possible by the freeness of $\Gamma$.
In fact, suppose that in a (necessarily split) extension $Z^{\times}\rightarrowtail G\twoheadrightarrow \Gamma$ of abelian groups, a section $g\mapsto\tilde{g}$, $\Gamma\rightarrow G$ affords a $2$-coboundary $\tau : \Gamma\times\Gamma\rightarrow Z^{\times}$.
The group homomorphism $f : Z^{\times}\times\Gamma\rightarrow G$ uniquely determined by
\[ f|_{Z^{\times}}=\mathrm{id},\quad f(g_{i})=\eta(g_{i})\tilde{g}_{i}\ \ (i\in I_{+}) \]
gives an equivalence $Z^{\times}\times\Gamma\xrightarrow{\simeq} G$ of extensions.
This implies $\tau=\partial\eta$.
We can apply this result for $\tau$ to $\sigma'\sigma^{-1}$, since
$\sigma'\sigma^{-1}$, being a $2$-coboundary, is symmetric.
\end{proof}

\begin{remark}
G\"unther \cite[Satz 4.9, Lemma 4.11]{G1} determined the set $\ZCleft(U_{q};Z)$, assuming that the GCM $\mathbb{A}$ is of finite type.
Our parametrization by $(u,\mu)$ is essentially the same as G\"unther's.
But, he supposes, with our notation, that $\mu$ can take non-zero values outside $\Xi(u)$, unless $Z$ is an integral domain; this is wrong.
\end{remark}

\begin{theorem}
Assume ($\mathrm{C}_{q}$), and that the characteristic $\ch k\neq 2$.
Then for any trivial $U_{q}'$-module $M$,
\[ H^{2}(U_{q}',M)=0. \]
\end{theorem}
\begin{proof}
Since the assumptions (C1)--(C4) are satisfied (see the proof of Theorem 7.14), it suffices by Theorem 6.4 to prove that $\mathcal{H}_{\varepsilon}^{\lambda}(\Gamma,\Xi;M)$ contains only one element.

We see from Example 7.11 that under ($C_{q}$), 
\[ \Xi = \Xi_{\lambda}=\{(i,-i)\;|\; i\in I_{+}\}. \]
Let $(s,m)\in\mathcal{Z}_{\varepsilon}^{\lambda}(\Gamma,\Xi;M)$.
Then, $s\in Z^{2}(\Gamma,M)$, $m=(m_{i,-i})_{i\in I_{+}}$, such that
\begin{equation}
2s(g,g_{i})=2s(g_{i},g)\quad (i\in I_{+},\ g\in\Gamma).
\end{equation}
Take another $(s',m')$.
We have to find a $1$-cochain $t : \Gamma\rightarrow M$ such that
\begin{eqnarray}
s' & = & s+\partial t, \\
m_{i,-i}' & = & -\lambda_{i,-i}2t(g_{i})+m_{i,-i} \quad (i\in I_{+}).
\end{eqnarray}
We can choose elements $t(g_{i})\in M$ ($i\in I_{+}$) satisfying (7.17), since we assume $\ch k\neq 2$.
The assumption and (7.15) imply that $s(g_{j},g_{i})=s(g_{i},g_{j})$ for all $i<j$, whence $s$ is symmetric (and a $2$-coboundary).
Just as in the last part of the proof of Theorem 7.14, $g_{i}\mapsto t(g_{i})$ can extend to a $1$-cochain $t$ satisfying (7.16).
\end{proof}

\begin{remark}
(1) Suppose $\ch k\neq 2$.
For any trivial $U_{q}'$-module $M$,
\[ H^{1}(U_{q}',M)=0 \]
since we see that an augmented algebra map $U_{q}'\rightarrow k_{M}$ is necessarily trivial; see \cite[Sect.\ 9.2]{W}.

(2) Suppose that the GCM $\mathbb{A}$ is of finite type.
We choose $\mathbb{D}=\diag(\cdots d_{i}\cdots)$ and $\mathfrak{h}_{\mathbb{Z}}$ in the standard way, so that $d_{i}\in\{1,2,3\}$, and $\mathfrak{h}_{\mathbb{Z}}=\bigoplus_{i=1}^{n}\mathbb{Z}d_{i}\alpha_{i}^{\vee}$. Hence, $U_{q}'=U_{q}$.
Suppose $\ch k =0$, and that $q$ is transcendental over $\mathbb{Q}$.
Then every finite-dimensional $U_{q}$-module is semisimple \cite[5.17 Theorem]{Ja}.
Let $M$ be a finite-dimensional simple $U_{q}$-module.
Then we see just as in the classical, Lie-algebra situation \cite[Theorem 7.8.9]{W}, that if $M\neq k$, then
\begin{equation}
H^{i}(U_{q},M)=\Ext_{U_{q}}^{i}(k,M)=0\quad (i\geq 0),
\end{equation}
since the center $Z(U_{q})$ of $U_{q}$ acts on $M$ through a non-trivial algebra map $Z(U_{q})\rightarrow k$, while it acts on $k$ through the trivial $Z(U_{q})\rightarrow k$; see \cite[6.26 Claim]{Ja}.
By Theorem 7.16 and Part 1 above, the equation (7.18) holds true when $i=1,2$, even if $M=k$. This result is a quantum analogue of Whitehead's lemmas \cite[Corollary XVIII.3.2]{K}, \cite[Corollaries 7.8.10, 7.8.12]{W}.

(3) It is possible to apply some classification result of non-central cleft extensions to compute the 2nd cohomology with non-trivial coefficients.
For example, let $T$ be Taft's Hopf algebra of dimension $N^{2}$, supposing that $k$ contains a primitive $N$-th root $\zeta$ of $1$.
If we let $\Gamma=\langle g\rangle$ denote the cyclic group of order $N$ with a generator $g$, and let $\chi\in\hat{\Gamma}$ such that $\chi(g)=\zeta$, then $T=k[x]/(x^{N})\dotrtimes\Gamma$ in which $x$ spans an object $\{x,g,\chi\}$ in ${}^{\Gamma}_{\Gamma}\YD$.
There are precisely $N$ simple $T$-modules $k(i)$ ($0\leq i <N$), which is $1$-dimensional and arises from the character $\chi^{i}$.
From the results in \cite{M1} on $T$-cleft extensions, we see
\[ H^{2}(T,k(i))=\left\{\begin{array}{ll} k & (i=0) \\ 0 & (0<i<N). \end{array}\right. \]
\c{S}tefan \cite[Theorem 2.8]{St} proved this result, supposing that $N$ is a prime.
\end{remark}

\section*{Appendix. On the Hopf algebras $u(\mathcal{D},\lambda,\mu)$ defined by Andruskiewitsch and Schneider}

\setcounter{equation}{0}
\renewcommand{\theequation}{A.\arabic{equation}}
\setcounter{theorem}{0}
\renewcommand{\thetheorem}{A.\arabic{theorem}}

In this appendix, we suppose that $k$ is a field of characteristic zero.
Very recently, Grunenfelder and Mastnak \cite[Theorem 3.5]{GM} formulated essentially the following.

\begin{theorem}
The finite-dimensional pointed Hopf algebra $u(\mathcal{D},\lambda,\mu)$, which was defined by Andruskiewitsch and Schneider \cite[Definition 4.3]{AS3} as a general form in their classification list, is a cocycle deformation of the associated graded Hopf algebra $u(\mathcal{D},0,0)$.
\end{theorem}

The proof by Grunenfelder and Mastnak depends on the author's \cite[Theorem 2]{M3}, which, however, was wrong, missing an assumption; the corrected formulation, given in the author's private notes, was published in \cite{BDR} as Theorem 3.4.
An earlier version of \cite{GM} missed the missing assumption, which was soon verified in this appendix in an earlier version of our paper.
The revised paper \cite{GM} refers to the appendix for verifying the assumption.
On this occasion of revision we will give a more self-contained (and hopefully complete) proof of Theorem A.1; see Remark A.8 below.

Let $\mathbb{A}=(a_{ij})_{i,j}$ be a $\theta\times\theta$ GCM of finite type, and let $\Phi^{+}$ denote the set of the positive roots in the associated root system.
Set $I=\{1,2,\dotsc,\theta\}$, and let $I_{1},\dotsc,I_{t}$ ($\subset I$) denote the connected components in the Dynkin diagram of $\mathbb{A}$.
Regard $I\subset\Phi^{+}$ as the subset of the simple roots.
Let $\Gamma$ be a finite abelian group, and suppose that we are given
\[ V=\{x_{i},g_{i},\chi_{i}\}_{i\in I}\in {}^{\Gamma}_{\Gamma}\YD. \]
Suppose that the matrix $q=(q_{ij})$ given by $q_{ij}=\chi_{j}(g_{i})$ is of Cartan type (see (6.1)) with respect to $\mathbb{A}$.
Suppose $q_{ii}\neq 1$ ($i\in I$), and that each $q_{ii}$ is a root of $1$ of odd order, satisfying the additional assumption given by \cite[(2.5)]{AS3}; it ensures that the order of $q_{ii}$ is constant on each component $I_{r}$.
Let $\alpha\in\Phi^{+}$, and let $I_{r(\alpha)}$ denote the connected component in which $\alpha$ has its support.
Let $N_{\alpha}$ denote the constant value, mentioned above, on $I_{r(\alpha)}$.
Andruskiewitsch and Schneider \cite[Sect.\ 2]{AS3} defined the root vector $x_{\alpha}$ ($\in TV$) as an iterated braided commutator of those $x_{i}$ with $i\in I_{r(\alpha)}$.
Let $\chi_{\alpha}\in\hat{\Gamma}$ defined by $gx_{\alpha}=\chi_{\alpha}(g)x_{\alpha}$ ($g\in\Gamma$).

To reconstruct the relevant algebras by our method, we first verify that $I=I_{1}\sqcup\dotsb\sqcup I_{t}$ is a partition of $I$ into mutually disconnected (in our sense), non-empty subsets.
We may suppose that $i>j$ whenever $r(i)>r(j)$ (or $|i|>|j|$ with our notation).
Then the set $\Theta$ turns to be
\[ \Theta=\{(i,j)\in I\times I\;|\;i>j,\ i\not\sim j\}, \]
where $i\sim j$ (resp., $i\not\sim j$) means that $i$ and $j$ are connected (resp., disconnected).
To construct the graded Hopf algebra $H^{0}$ from the $V$ above, choose $L\in\mathbb{N}$ so that $L>\mathrm{ht}(\alpha)$ for all $\alpha\in\Phi^{+}$, and define
\[ \emptyset=\mathcal{R}_{(0)}\subset\mathcal{R}_{(1)}\subset\dotsb\subset\mathcal{R}_{(L)}=\mathcal{R} \]
by
\begin{eqnarray*}
\mathcal{R}_{(1)} & = & \{(\ad_{c}x_{i})^{1-a_{ij}}(x_{j})\;|\;i\sim j,\ i\neq j\}, \\
\mathcal{R}_{(h)} & = & \mathcal{R}_{(1)}\cup\{x_{\alpha}^{N_{\alpha}}\;|\;\alpha\in\Phi^{+},\ \mathrm{ht}(\alpha)<h,\ \chi_{\alpha}^{N_{\alpha}}\neq 1\}\quad (1<h\leq L).
\end{eqnarray*}
By \cite[Lemma A.1]{AS1}, $\mathcal{R}_{(1)}$ consists of primitives in $TV$.

Let $\alpha\in\Phi^{+}$.
Let $\beta_{1},\dotsc,\beta_{p}$ be the elements in $\Phi^{+}$ with support in $I_{r(\alpha)}$ which are given in such a special order as described in \cite[Sect.\ 2]{AS3}.
By \cite[(2.12)]{AS3}, the coproduct $\Delta(x_{\alpha}^{N_{\alpha}})$ of $TV/(\mathcal{R}_{(1)})$ is given by
\begin{equation}
\Delta(x_{\alpha}^{N_{\alpha}})=x_{\alpha}^{N_{\alpha}}\otimes 1+1\otimes x_{\alpha}^{N_{\alpha}}+\sum_{b,c}t_{b,c}z^{b}\otimes z^{c},
\end{equation}
in which the summation is taken over all non-zero vectors
\[ b=(b_{1},\dotsc,b_{p}),\ \ c=(c_{1},\dotsc,c_{p})\ \ \mbox{in}\ \ \mathbb{N}^{p}\]
such that $\sum_{l=1}^{p}(b_{l}+c_{l})\beta_{l}=\alpha$.
In addition, $t_{b,c}\in k$, and for $a=(a_{1},\dotsc,a_{p})\in\mathbb{N}^{p}$,
\begin{equation}
z^{a}:=x_{\beta_{1}}^{Na_{1}}\dotsm x_{\beta_{p}}^{Na_{p}},\ \ \mbox{where}\ \ N:=N_{\alpha}.
\end{equation}
In (A.1), the heights $\mathrm{ht}(\beta_{l})$ of those $\beta_{l}$ for which $x_{\beta_{l}}^{N}$ appears in $z^{b}$ or $z^{c}$ are all less than $\mathrm{ht}(\alpha)$.
By using \cite[Theorem 2.6.1]{AS3}, one sees that if $\chi_{\alpha}^{N_{\alpha}}\neq 1$, or namely if $x_{\alpha}^{N_{\alpha}}$ is not $\Gamma$-invariant, some $x_{\beta_{l}}^{N}$ that appears in $z^{b}$ or $z^{c}$ is not $\Gamma$-invariant.
Therefore every element in $\mathcal{R}_{(h)}$ ($h>1$) is a primitive modulo $(\mathcal{R}_{(h-1)})$.
We can thus define as in Definition 2.2, a graded Hopf algebra,
\[ H^{0}=TV/(\mathcal{R}^{0})\dotrtimes\Gamma. \]

Choose arbitrarily $\lambda_{ij}\in k$ for each $(i,j)\in\Theta$ with $\chi_{i}\chi_{j}=1$, and set $\lambda=(\lambda_{ij})\in k^{\Xi}$.
By Proposition 4.2.3, we obtain a Hopf algebra,
\[ H^{\lambda}=TV\dotrtimes\Gamma/(\mathcal{R}^{\lambda}), \]
and an $(H^{0},H^{\lambda})$-bicleft object,
\[ A(\lambda)=TV\rtimes\Gamma/(\mathcal{R}(\lambda)). \]

\begin{remark}
We can define $\mathcal{R}_{(h)}$ ($1<h\leq L$) alternatively so as
\[ \mathcal{R}_{(h)}=\mathcal{R}_{(1)}\cup\{x_{\alpha}^{N_{\alpha}}\;|\;\alpha\in\Phi^{+},\ \mathrm{ht}(\alpha)<h\}. \]
In this case, our Theorem 4.3, telling that $H^{\lambda}$ and $H^{0}$ are cocycle deformations of each other, is precisely Didt's \cite[Theorem 1]{Di}.
\end{remark}

Recall from Remark 2.6 and an analogous fact that $H^{0}$, $H^{\lambda}$ and $A(\lambda)$ all naturally include the Hopf algebras
\begin{equation}
E_{r}=T_{r}\dotrtimes\Gamma\ \ (1\leq r\leq t),\quad \mbox{where}\ \ T_{r}=\frac{TV_{r}}{(\mathcal{R}_{r})},
\end{equation}
so that the coactions
\begin{equation}
H^{0}\otimes A(\lambda)\xleftarrow{\rho}A(\lambda)\xrightarrow{\rho'}A(\lambda)\otimes H^{\lambda}
\end{equation}
restricts to the coproduct of $E_{r}$.

Define
\[ \Phi_{0}^{+} := \{\alpha\in\Phi^{+}\;|\;\chi_{\alpha}^{N_{\alpha}}=1\}. \]
Suppose $\alpha\in\Phi_{0}^{+}$.
Let $g_{\alpha}=|x_{\alpha}|$ ($\in\Gamma$), the degree of $x_{\alpha}$ in $\Gamma$, and set $N=N_{\alpha}$.
It follows by (A.1) that the coproducts $\Delta(x_{\alpha}^{N})$ both on $H^{0}$ and on $H^{\lambda}$, $\rho(x_{\alpha}^{N})$ and $\rho'(x_{\alpha}^{N})$ are all given by the same formula
\begin{equation}
x_{\alpha}^{N}\otimes 1+g_{\alpha}^{N}\otimes x_{\alpha}^{N}+\sum_{b,c}t_{b,c}z^{b}h^{c}\otimes z^{c},
\end{equation}
where
\[ h^{c}:=g_{\beta_{1}}^{Nc_{1}}\dotsm g_{\beta_{p}}^{Nc_{p}}, \]
and $b,c$ now may be restricted to those for which $b_{l}=c_{l}=0$ unless $\beta_{l}\in\Phi_{0}^{+}$.

\begin{lemma}
Let $\alpha\in\Phi_{0}^{+}$.
Then, $x_{\alpha}^{N_{\alpha}}$ is central all in $H^{0}$, in $A(\lambda)$ and in $H^{\lambda}$.
\end{lemma}
\begin{proof}
The result in $H^{0}$ follows by \cite[Theorem 2.6.3]{AS3}.

To prove the result in $A(\lambda)$, set $N=N_{\alpha}$.
Obviously, $gx_{\alpha}^{N}=x_{\alpha}^{N}g$ for all $g\in\Gamma$.
Choose any $i\in I$.
We will prove $x_{i}x_{\alpha}^{N}=x_{\alpha}^{N}x_{i}$.
This follows from the result in $H^{0}$ and Remark 2.6, if $i\in I_{r(\alpha)}$.
We suppose $i\not\in I_{r(\alpha)}$, to prove the desired equation by induction on $\mathrm{ht}(\alpha)$.
Since $x_{\alpha}^{N}$ is central in $H^{0}$, $\rho(x_{i})$ ($=x_{i}\otimes 1+g_{i}\otimes x_{i}$) commutes with the first term $x_{\alpha}^{N}\otimes 1$ of $\rho(x_{\alpha})$ given in (A.5).
It commutes with the last summation, by the induction hypothesis.
Here, notice from (6.1) that in $H^{0}$ (and also in $A(\lambda)$, $H^{\lambda}$),
\begin{equation}
x_{i}h^{c}=\chi_{\beta_{1}}^{Nc_{1}}\dotsm\chi_{\beta_{p}}^{Nc_{p}}(g_{i})=h^{c}x_{i}.
\end{equation}
If we set $w=x_{i}x_{\alpha}^{N}-x_{\alpha}^{N}x_{i}$ ($\in A(\lambda)$), it follows that $\rho(w)=g_{i}g_{\alpha}^{N}\otimes w$, whence $w=dg_{i}g_{\alpha}^{N}$ for some $d\in k$.
But, since $N>1$, $d$ must be zero, as is seen from (2.18).
Hence, $x_{i}x_{\alpha}^{N}=x_{\alpha}^{N}x_{i}$, which completes the induction.

To prove the result in $H^{\lambda}$, we will modify the proof in $A(\lambda)$ just above, and prove by induction on $\mathrm{ht}(\alpha)$ that for $i\not\in I_{r(\alpha)}$, $w=x_{i}x_{\alpha}^{N}-x_{\alpha}^{N}x_{i}$ ($\in H^{\lambda}$) equals zero.
By using the result in $A(\lambda)$ and the induction hypothesis, we see $\rho'(w)=g_{i}g_{\alpha}^{N}\otimes w$, whence $w=\varepsilon(w)g_{i}g_{\alpha}^{N}=0$, as desired.
\end{proof}

Let $B$ denote the subalgebra of $H^{0}$ generated by all $x_{\alpha}^{N_{\alpha}}$ ($\alpha\in\Phi_{0}^{+}$);
this is a left coideal subalgebra.
Let $\Gamma'$ denote the subgroup of $\Gamma$ generated by all $g_{\alpha}^{N_{\alpha}}$ ($\alpha\in\Phi_{0}^{+}$).
Notice that $B$ and $\Gamma'$ generate a Hopf subalgebra of $H^{0}$.

\begin{proposition}
(1) $B$ is a commutative polynomial algebra with indeterminates $x_{\alpha}^{N_{\alpha}}$ ($\alpha\in\Phi_{0}^{+}$).

(2) $H^{0}$ and $H^{\lambda}$ both naturally include $B\otimes k\Gamma'$ as a central Hopf subalgebra.

(3) $A(\lambda)$ naturally includes $B\otimes k\Gamma'$ as a central subalgebra, so that the coactions $\rho,\rho'$ given in (A.4) restrict to the coproduct on $B\otimes k\Gamma'$.
\end{proposition}
\begin{proof}
(1) We see from \cite[Theorem 2.6]{AS3} that those $x_{\alpha}^{N_{\alpha}}$ ($\alpha\in\Phi_{0}^{+}$) which are contained in a single $T_{r}$ are algebraically independent.
This together with Lemma A.3 and Proposition 2.5.1 imply the desired result.

(2), (3) First, notice from (A.6) that each $g_{\alpha}^{N_{\alpha}}$ ($\alpha\in\Phi_{0}^{+}$) is central all in $H^{0}$, $A(\lambda)$ and $H^{\lambda}$.
Then the desired results follow from Lemma A.3, Propositions 2.5.1 and 4.2.2; see also Remark 2.6.
\end{proof}

Choose arbitrarily an algebra map $f : B\rightarrow k$, and extend it uniquely to an algebra map
\[ f : B\otimes k\Gamma'\rightarrow k\ \ \mbox{with}\ \ f(g)=1\ \ (g\in\Gamma'), \]
which is necessarily convolution-invertible.
One sees that
\[ f\cdot h=\sum h_{1}f(h_{2}),\quad h\cdot f^{-1}=\sum f^{-1}(h_{1})h_{2}\quad (h\in B\otimes k\Gamma') \]
define algebra automorphisms of $B\otimes k\Gamma'$.
Notice that
\begin{equation}
h\mapsto f\cdot h\cdot f^{-1},\quad B\otimes k\Gamma'\xrightarrow{\simeq}B\otimes k\Gamma'
\end{equation}
is a Hopf algebra automorphism.
Define three ideals of $B\otimes k\Gamma'$ by
\[ \mathfrak{a}:=B^{+}\otimes k\Gamma',\quad \mathfrak{b}:=f\cdot\mathfrak{a},\quad \mathfrak{c}:=f\cdot\mathfrak{a}\cdot f^{-1}, \]
and let
\[ (\mathfrak{a})\subset H^{0},\quad (\mathfrak{b})\subset A(\lambda),\quad (\mathfrak{c})\subset H^{\lambda} \]
denote the ideals generated by $\mathfrak{a},\mathfrak{b},\mathfrak{c}$ in the respective algebras.
Notice that $\mathfrak{a}$, $\mathfrak{c}$, $(\mathfrak{a})$ and $(\mathfrak{c})$ are all Hopf ideals.
Let $\varphi : B\otimes k\Gamma'\rightarrow k\Gamma'$ be a unique algebra map that makes the following commute; cf.\ $\varphi_{\mu}$ defined in \cite[Theorem 2.13]{AS3}.
\begin{center}
\begin{picture}(120,50)
\put(0,35){$B\otimes k\Gamma'$}
\put(40,37){\vector(1,0){30}}
\put(75,35){$B\otimes k\Gamma'$}
\put(43,43){(A.7)}
\put(50,28){$\simeq$}
\put(30,30){\vector(1,-1){20}}
\put(80,30){\vector(-1,-1){20}}
\put(11,17){$\varepsilon\otimes\mathrm{id}$}
\put(74,17){$\varphi$}
\put(48,0){$k\Gamma'$}
\end{picture}
\end{center}
Here, $\varepsilon\otimes\mathrm{id}$, arising from the braided Hopf algebra map $\varepsilon : B\rightarrow k$, is a Hopf algebra map.
Hence, $\varphi$ is a Hopf algebra map, which is identical on $k\Gamma'$.
It is easy to prove the following lemma.

\begin{lemma}
$\mathfrak{c}=\Ker\varphi$.
This ideal is generated by $x_{\alpha}^{N_{\alpha}}-\varphi(x_{\alpha}^{N_{\alpha}})$ ($\alpha\in\Phi_{0}^{+}$).
\end{lemma}

\begin{proposition}
The coactions $\rho,\rho'$ in (A.4) induce
\[ H^{0}/(\mathfrak{a})\otimes A(\lambda)/(\mathfrak{b})\leftarrow A(\lambda)/(\mathfrak{b})\rightarrow A(\lambda)/(\mathfrak{b})\otimes H^{\lambda}/(\mathfrak{c}), \]
by which $A(\lambda)/(\mathfrak{b})$ is an $(H^{0}/(\mathfrak{a}),H^{\lambda}/(\mathfrak{c}))$-bicleft object.
Hence, $H^{0}/(\mathfrak{a})$ and $H^{\lambda}/(\mathfrak{c})$ are cocycle deformations of each other.
\end{proposition}
\begin{proof}
It suffices to prove that $A(\lambda)/(\mathfrak{b})$ is a biGalois object (i.e., a biGalois extension over $k$), which is necessarily a bicleft object since $H^{0}/(\mathfrak{a})$ and $H^{\lambda}/(\mathfrak{c})$ are pointed.

We know that $f\cdot : B\rightarrow B\otimes k\Gamma'\subset A(\lambda)$ is an $H^{0}$-comodule algebra map from a central left coideal subalgebra $B$ of a pointed Hopf algebra $H^{0}$, which is faithfully flat over $B$, into the center of an $H^{0}$-Galois object $A(\lambda)$.
It follows by (the opposite side version of) G\"unther's theorem \cite[Theorem 4]{G2} that $A(\lambda)/(f\cdot B^{+})=A(\lambda)/(\mathfrak{b})$ is a Galois object for $H^{0}/(B^{+})=H^{0}/(\mathfrak{a})$ with respect to the coaction induced from $\rho$.

To see that $A(\lambda)/(\mathfrak{b})$ is a right $H^{\lambda}/(\mathfrak{c})$-Galois object, we will prove that the (so-called Galois) isomorphism
\[ \gamma : A(\lambda)\otimes A(\lambda)\xrightarrow{\simeq} A(\lambda)\otimes H^{\lambda} \]
given by $\gamma(a\otimes b)=a\rho'(b)$ maps
\[ X:=(\mathfrak{b})\otimes A(\lambda)+A(\lambda)\otimes (\mathfrak{b}) \]
onto
\[ Y:=(\mathfrak{b})\otimes H^{\lambda}+A(\lambda)\otimes (\mathfrak{c}). \]
Set $B'=B\otimes k\Gamma'$.
Since we see as in the proof of \cite[Theorem 2]{M3} that
\begin{equation}
\gamma(\mathfrak{b}\otimes B'+B'\otimes\mathfrak{b})=\mathfrak{b}\otimes B'+B'\otimes\mathfrak{c},
\end{equation}
it easily follows that $\gamma(X)\subset Y$.
Let
\[ \phi=\phi^{\lambda} : H^{\lambda}\rightarrow A(\lambda) \]
denote the section as given in Proposition 4.2.2.
The composite inverse $\gamma^{-1}$ of $\gamma$ is given by
\[ \gamma^{-1}(a\otimes h)=\sum a\phi^{-1}(h_{1})\otimes\phi(h_{2})\quad (a\in A(\lambda),\ \ h\in H^{\lambda}). \]
We see then easily that $\gamma^{-1}((\mathfrak{b})\otimes H^{\lambda})\subset (\mathfrak{b})\otimes A(\lambda)\subset Y$.
By Lemma A.5, it remains to prove that for $\alpha\in\Phi_{0}^{+}$,
\[ \gamma^{-1}(1\otimes (x_{\alpha}^{N_{\alpha}}-\varphi(x_{\alpha}^{N_{\alpha}}))H^{0})\subset Y. \]
But, this follows from (A.8) and the next lemma.
\end{proof}

\begin{lemma}
Let $\phi : H^{\lambda}\rightarrow A(\lambda)$ be as above.
For each $1\leq r\leq t$, let $h_{r}$ be an arbitrary element in the Hopf algebra $E_{r}$ as given in (A.3), which is naturally embedded both into $A(\lambda)$, and into $H^{\lambda}$.
\begin{enumerate}
\renewcommand{\labelenumi}{(\arabic{enumi})}
\item $\phi(h_{1}\dotsm h_{t})=h_{1}\dotsm h_{t}$.
\item $\phi^{-1}(h_{1}\dotsm h_{t})=S(h_{t})\dotsm S(h_{1})$, where $S$ denotes the antipodes of $E_{r}$.
\end{enumerate}
\end{lemma}
\begin{proof}
Recall that $H^{\lambda}$ and $A(\lambda)$ are both naturally identified with $T_{1}\otimes\dotsb\otimes T_{t}\otimes k\Gamma$, where $T_{r}$ are as in (A.3).
Let $y_{r}\in T_{r}$ and $g\in\Gamma$.
Then, $\phi$ is given by
\[ \phi(y_{1}\dotsm y_{t}g)=y_{1}\dotsm y_{t}g. \]

(1) We can prove by downward induction on $1\leq r\leq t$ that $\phi(h_{r}\dotsm h_{t})=h_{r}\dotsm h_{t}$.

(2) Define $\psi : H^{\lambda}\rightarrow A(\lambda)$ by
\[ \psi(y_{1}\dotsm y_{t}g)=g^{-1}S(y_{t})\dotsm S(y_{1}). \]
Then one sees just as above that
\[ \psi(h_{1}\dotsm h_{t})=S(h_{t})\dotsm S(h_{1}), \]
which implies $\psi=\phi^{-1}$.
This proves the desired result.
\end{proof}

\begin{proof}[Proof of Theorem A.1]
By \cite[Definition 4.3]{AS3}, the Hopf algebras $u(\mathcal{D},\lambda,\mu)$ are parameterized by
\[ \lambda=(\lambda_{ij})\in k^{\Xi},\quad \mu=(\mu_{\alpha})\in k^{\Phi_{0}^{+}} \]
in which $\mu_{\alpha}=0$ if $g_{\alpha}^{N_{\alpha}}=1$.
The graded Hopf algebra $u(\mathcal{D},0,0)$ is given by the zero parameters $\lambda_{ij}=\mu_{\alpha}=0$.
We see easily that
\[ u(\mathcal{D},0,0)=H^{0}/(\mathfrak{a}). \]

Choose arbitrarily a connected component $I_{r}$ in $I$, and let $\beta_{1},\dotsc,\beta_{p}$ be, as above, the positive roots with support in $I_{r}$.
For $0\neq a\in\mathbb{N}^{p}$, let $z^{a}$ be as in (A.2), which, however is now regarded as an element in $B$ (therefore, we may suppose $a_{l}=0$ if $\chi_{\beta_{l}}^{N}\neq 1$).
For the algebra map $f : B\otimes k\Gamma'\rightarrow k$ as above, define an element in $B\otimes k\Gamma'$ by
\[ \delta(z^{a}):=z^{a}+f^{-1}(z^{a})(1-h^{a})+\sum_{b,c}t_{b,c}^{a}f^{-1}(z^{b})z^{c}, \]
where $t_{b,c}^{a}$ is the scalar as defined in \cite[Lemma 2.8]{AS3}, and the summation $\sum_{b,c}$ is taken over all non-zero vectors $b,c$ in $\mathbb{N}^{p}$ such that $\sum_{l=1}^{p}(b_{l}+c_{l})\beta_{l}=\sum_{l=1}^{p}a_{l}\beta_{l}$.
As was essentially proved by \cite[Proposition 3.3.6]{GM}, we have
\[ f\cdot z^{a}\cdot f^{-1}=\delta(z^{a})+\sum_{b,c}t_{b,c}^{a}\delta(z^{b})h^{c}f(z^{c}). \]
By definition of the Hopf algebra map $\varphi : B\otimes k\Gamma'\rightarrow k\Gamma'$, we have $\varphi(f\cdot z^{a}\cdot f^{-1})=0$.
By induction on the height $\mathrm{ht}(\sum_{l=1}^{p}a_{l}\beta_{l})$, we see $\varphi(\delta(z^{a}))=0$, or namely
\begin{equation}
\varphi(z^{a})+f^{-1}(z^{a})(1-h^{a})+\sum_{b,c}t_{b,c}^{a}f^{-1}(z^{b})\varphi(z^{c})=0.
\end{equation}
(We can also conclude that if $h^{a}=1$, then the value $f^{-1}(z^{a})$ does not have any effect on the Hopf algebra automorphism (A.7).)
Define $f$ so that
\[ f^{-1}(x_{\alpha}^{N_{\alpha}})=-\mu_{\alpha}\quad (\alpha\in\Phi_{0}^{+}). \]
(This is possible, since $f^{-1}=f\circ S$, where $S$ is the antipode, that is an involution, of $B\otimes k\Gamma'$.)
Compare (A.9) with \cite[(2.15)]{AS3}, setting $f^{-1}(z^{a})=-\mu_{a}$.
Then we see from \cite[Theorem 2.13.1]{AS3} that $\varphi(x_{\alpha}^{N_{\alpha}})$ coincides with the element $u_{\alpha}(\mu)$ defined by \cite[Definition 2.14]{AS3}.
Moreover, in the definition of $H^{\lambda}$, replace $\lambda$ with $-\lambda=(-\lambda_{ij})$.
Then we see
\[ u(\mathcal{D},\lambda,\mu)=H^{-\lambda}/(\mathfrak{c}). \]
The desired result now follows by Proposition A.6.
\end{proof}

\begin{remark}
We have modified (and supplemented) the proof of Grunenfelder and Mastnak \cite[Theorem 3.5]{GM}.
Their proof starts from an algebra map $f : K(\mathcal{D})\dotrtimes\Gamma\rightarrow k$, where $K(\mathcal{D})\dotrtimes\Gamma$ is a Hopf algebra which has our $B\otimes k\Gamma'$ as a sub-quotient; see the paragraph preceding \cite[Theorem 3.5]{GM}.
But, the well-definedness of their $f$ seems unclear (at least for the author).
They refer to (an earlier version of) our paper for the proof that the bicomodule algebra, our $A(\lambda)/(\mathfrak{b})$, which is to be a biGalois object is non-zero;
this is the assumption for \cite[Theorem 3.4]{BDR} which was originally missing in \cite{M3}.
\end{remark}

\end{document}